\documentclass[12pt]{amsart}

\usepackage{amsmath}
\usepackage{amssymb}
\usepackage{latexsym}
\usepackage{graphicx} 
\usepackage{mathrsfs}
\usepackage{mathtools}
\usepackage{amsthm}
\usepackage{tikz}
\usetikzlibrary{cd}
\allowdisplaybreaks[4]

\usepackage[colorlinks,citecolor=blue,linkcolor=blue,linktocpage,unicode]{hyperref}

\numberwithin{equation}{section}
\numberwithin{figure}{section}

\newtheorem{thm}{Theorem}[section]
\newtheorem{prop}[thm]{Proposition}
\newtheorem{lemma}[thm]{Lemma}
\newtheorem{cor}[thm]{Corollary}

\theoremstyle{definition}
\newtheorem{example}[thm]{Example}
\newtheorem{defin}[thm]{Definition}
\newtheorem{remark}[thm]{Remark}
\newtheorem*{ack}{Acknowledgements}

\makeatletter
\renewenvironment{proof}[1][\proofname]{\par
  \pushQED{\qed}%
  \normalfont \topsep6\p@\@plus6\p@\relax
  \trivlist
  \item\relax
  {
  #1\@addpunct{.}}\hspace\labelsep\ignorespaces
}{%
  \popQED\endtrivlist\@endpefalse
}
\makeatother

\title[Torus-equivariantly embedded toric manifolds]{Torus-equivariantly embedded toric manifolds associated to affine subspaces}
\author{Kentaro Yamaguchi}

\address{Department of Mathematical Sciences, Tokyo Metropolitan University, 1-1 Minami-Ohsawa, Hachioji, Tokyo, 192-0397, Japan}
\email{yamaguchi-kentaro@ed.tmu.ac.jp}

\begin{document}

\begin{abstract}
We study the closure of a complex subtorus in a toric manifold.
If the closure of the complex subtorus is a smooth complex submanifold in the toric manifold, then the subtorus action on such submanifold is Hamiltonian.
In this case, we may think of the embedding of the submanifold as torus-equivariant. 
We show that the image of the moment map for the Hamiltonian subtorus action on our submanifold coincides with the image of the Delzant polytope of the ambient toric manifold under the pullback of the inclusion of the tori.
The submanifolds constructed in the present paper are called \textit{torus-equivariantly embedded toric manifolds} with respect to the subtorus action.
\end{abstract}

\maketitle
%\tableofcontents

\section{Introduction}
\label{sec: introduction}

Delzant \cite{MR984900} established a one-to-one correspondence between compact symplectic toric manifolds and certain convex polytopes known as \textit{Delzant polytopes}.
Given a $2n$-dimensional compact symplectic toric manifold $X$, the image of a moment map for the Hamiltonian $T^{n}$-action on $X$ is a Delzant polytope $\Delta$ in $(\mathfrak{t}^{n})^{\ast} \cong \mathbb{R}^{n}$.
Conversely, given a Delzant polytope $\Delta$ in $(\mathfrak{t}^{n})^{\ast}$, we can construct a compact symplectic toric manifold $X$ whose moment polytope is $\Delta$. 
This construction is called the \textit{Delzant construction}.
From the Delzant construction, symplectic toric manifolds are canonically equipped with a K\"{a}hler structure \cite{MR1293656,MR2039163}.
We can identify the complements of toric divisors in a symplectic toric manifold $X$ with a complex torus $(\mathbb{C}^{\ast})^{n}$, whose description allows us to consider complex coordinates in $X$.

\subsection{Main Results}
In this paper, we study complex submanifolds in compact toric manifolds $X$.
From a $k$-dimensional affine subspace $V$ in $\mathfrak{t}^{n} \cong \mathbb{R}^{n}$, we first construct a $k$-dimensional complex submanifold $C(V)$ in the toric divisor complements $\check{M} \cong (\mathbb{C}^{\ast})^{n}$ of the toric manifold $X$. This construction is inspired by \cite{MR3856842}, and $C(V) \cong (\mathbb{C}^{\ast})^{k}$ as Yamamoto noted there. 
In fact, $C(V) \cong (\mathbb{C}^{\ast})^{k}$ is a complex subtorus of $\check{M} \cong (\mathbb{C}^{\ast})^{n}$.
We then consider the conditions of $V$ when the (Zariski) closure $\overline{C(V)}$ is a $k$-dimensional complex submanifold in the toric manifold $X$ (Section \ref{subsec: construction of toric submanifolds}).
While $C(V)$ is a complex submanifold in $\check{M} \cong (\mathbb{C}^{\ast})^{n}$ for arbitrary affine subspace $V$ as Yamamoto showed in \cite[Lemma 6.1]{MR3856842}, $\overline{C(V)}$ may not be a complex submanifold in $X$ (see Example \ref{eg: cp2, n.g.} and Section \ref{sec: example}).

Suppose that $\overline{C(V)}$ is a smooth complex submanifold in $X$.
We then discuss the nature of the submanifolds $\overline{C(V)}$.
Toric manifolds $X$ are naturally equipped with a moment map $\mu: X \to (\mathfrak{t}^{n})^{\ast}$ for the $T^{n}$-action on them.
We can define the injective group homomorphism $i_{V}: T^{k} \to T^{n}$ by the data of $V$ (see Equation \ref{eq: inclusion of torus}).
Because the $T^{n}$-action on $X$ and $i_{V}: T^{k} \to T^{n}$ induce the $T^{k}$-action on $\overline{C(V)}$, we can determine the moment map $\overline{\mu}:\overline{C(V)} \to (\mathfrak{t}^{k})^{\ast}$ by $\overline{\mu} = i_{V}^{\ast}\circ \mu \circ i$, where $i:\overline{C(V)} \to X$ is the embedding (see Section \ref{subsec: torus actions on toric submanifolds} for detail).
We obtain the following diagram:
\begin{equation*}
	\begin{tikzcd}
		X \arrow[r,"\mu"] & (\mathfrak{t}^{n})^{\ast} \arrow[d, "i_{V}^{\ast}"] \\
		\overline{C(V)} \arrow[u, hook, "i"] \arrow[r, "\overline{\mu}"] & (\mathfrak{t}^{k})^{\ast}.
	\end{tikzcd}	
\end{equation*}
We compare the image of the moment map $\overline{\mu}:\overline{C(V)} \to (\mathfrak{t}^{k})^{\ast}$ for the $T^{k}$-action on our complex submanifold $\overline{C(V)}$ with the image of the moment map $\mu: X \to (\mathfrak{t}^{n})^{\ast}$ for the $T^{n}$-action on the ambient toric manifold.
\begin{thm}[Theorem \ref{thm: moment polytope of toric submanifolds}]
	Let $i_{V}:T^{k} \to T^{n}$ be an injective group homomorphism determined by a given affine subspace $V$ in $\mathfrak{t}^{n} \cong \mathbb{R}^{n}$.
	Assume that $\overline{C(V)}$ is a complex submanifold in $X$.
	Then the image of $\overline{\mu}$ is equal to the image of $i_{V}^{\ast}\circ \mu$, i.e., $\overline{\mu}(\overline{C(V)}) = (i_{V}^{\ast}\circ \mu)(X)$.
\end{thm}
We call the complex submanifolds $\overline{C(V)}$ \textit{torus-equivariantly embedded toric manifolds}.

\subsection{Outline}
This paper is organized as follows. 
In Section \ref{sec: coordinates chart}, 
we construct a system of complex coordinate charts on complex manifolds from matrices in $SL(n;\mathbb{Z})$. This construction helps us to consider a system of complex coordinate charts in toric manifolds.
In Section \ref{sec: toric manifold}, 
we review Delzant construction and construct a system of the inhomogeneous coordinate charts on compact toric manifolds using the construction established in Section \ref{sec: coordinates chart}.
In Section \ref{sec: toric submanifolds}, 
we give the conditions where the closure $\overline{C(V)}$ of a complex subtorus $C(V)$ is a complex submanifold in the ambient toric manifold.
Moreover, we consider the moment maps for the subtorus action on our complex submanifolds and compare them with the moment maps for the torus action on ambient toric manifolds. 
In Section \ref{sec: example}, we demonstrate some examples of torus-equivariantly embedded toric manifolds.

Throughout the paper, we express vectors as column vectors.

\section{Construction of Coordinate Charts from Matrices in $SL(n;\mathbb{Z})$}
\label{sec: coordinates chart}

This section is about construction of a system of complex coordinate charts from matrices in the special linear group $SL(n; \mathbb{Z})$.
Our idea is similar to \textit{the coordinate transformations} for compact toric manifolds by Duistermaat and Pelayo \cite{MR2496353}, but here we construct a system of complex coordinate charts in a general situation.

Let $\Lambda$ be a set and $Q^{\lambda} \in SL(n;\mathbb{Z})$ a matrix corresponding to each $\lambda \in \Lambda$ and $\mathbb{C}^{n}_{\lambda} = \{ z^{\lambda} = (z^{\lambda}_{1}, \ldots, z^{\lambda}_{n})\in \mathbb{C}^{n} \}\cong \mathbb{C}^{n}$ for each $\lambda \in \Lambda$. 
We define a matrix $D^{\lambda\mu} = (Q^{\lambda})^{-1}Q^{\mu} (=[d^{\lambda\mu}_{ij}])$ for any $\lambda,\mu \in \Lambda$ 
and a subset $U_{\lambda\mu} \subset \mathbb{C}^{n}_{\lambda}$ by 
\begin{equation*}
	U_{\lambda\mu} = \{ z^{\lambda} \in \mathbb{C}_{\lambda}^{n} \mid z^{\lambda}_{j} \neq 0 ~\text{if}~ d^{\lambda\mu}_{jl} < 0 ~\text{for some $l = 1, \ldots, n$}\}.
\end{equation*}
We introduce an equivalence relation on $\{U_{\lambda\mu}\}_{\lambda,\mu \in \Lambda}$.

\begin{defin}
\label{def: equiv relation}
	For $z^{\lambda} \in U_{\lambda\mu} \subset \mathbb{C}^{n}_{\lambda}$ and $z^{\mu} \in U_{\mu\lambda} \subset \mathbb{C}^{n}_{\mu} $, 
	we define a binary relation $z^{\lambda} \sim z^{\mu}$ by
	\begin{equation*}
		(z^{\mu}_{1},\ldots, z^{\mu}_{n}) = \left( \prod_{j=1}^{n} \left(z^{\lambda}_{j}\right)^{d^{\lambda\mu}_{j1}}, \ldots, \prod_{j=1}^{n} \left(z^{\lambda}_{j}\right)^{d^{\lambda\mu}_{jn}} \right).
	\end{equation*}
\end{defin}

\begin{prop}
\label{prop: equiv relation}
	The binary relation $\sim$ defined in Definition \ref{def: equiv relation} is an equivalence relation.
	
\begin{proof}
	We check that the binary relation $\sim$ satisfies the definition of equivalence relations.
	
	Since we define $D^{\lambda \mu} = (Q^{\lambda})^{-1} Q^{\mu}$ for $\lambda,\mu \in \Lambda$, we get $D^{\lambda\lambda} = (Q^{\lambda})^{-1} Q^{\lambda} = E_{n}$, where $E_{n}$ is the identity matrix.
	Hence, $z^{\lambda} = z^{\lambda}$, which means that $z^{\lambda} \sim z^{\lambda}$.

	Since $D^{\lambda \mu} = (Q^{\lambda})^{-1} Q^{\mu}$, we have $D^{\mu \lambda} = (Q^{\mu})^{-1} Q^{\lambda} = (D^{\lambda\mu})^{-1}$.
	Suppose $z^{\lambda} \sim z^{\mu}$, then we see $z^{\mu}_{i} = \prod_{j=1}^{n}(z^{\lambda}_{j})^{d^{\lambda\mu}_{ji}}$ for $i= 1, \ldots, n$.
	For $k = 1, \ldots, n$, we have 
		\begin{align*}
			\prod_{i=1}^{n} \left(z^{\mu}_{i}\right)^{d^{\mu\lambda}_{ik}}
			= \prod_{i=1}^{n} \prod_{j=1}^{n} \left( z^{\lambda}_{j} \right)^{d^{\lambda\mu}_{ji} d^{\mu\lambda}_{ik}} 
			= \prod_{j=1}^{n} \left( z^{\lambda}_{j} \right)^{\delta_{jk}} 
			= z^{\lambda}_{k}.
		\end{align*}
		Hence we obtain $z^{\mu} \sim z^{\lambda}$.
			
		For $\lambda,\mu, \sigma \in \Lambda$, we have
		\begin{equation*}
			D^{\lambda\sigma} = (Q^{\lambda})^{-1} Q^{\sigma} = (Q^{\lambda})^{-1} Q^{\mu} (Q^{\mu})^{-1} Q^{\sigma} = D^{\lambda\mu} D^{\mu\sigma}.
		\end{equation*}
		Suppose $z^{\lambda} \sim z^{\mu}, z^{\mu} \sim z^{\sigma}$, then for $i=1,\ldots,n$ we have
		\begin{align*}
			z^{\sigma}_{i} 
			= \prod_{j=1}^{n} \left( z^{\mu}_{j} \right)^{d^{\mu\sigma}_{ji}} 
			= \prod_{j=1}^{n} \left( \prod_{k=1}^{n} \left( z^{\lambda}_{k} \right)^{d^{\lambda\mu}_{kj}} \right)^{d^{\mu\sigma}_{ji}} 
			= \prod_{k=1}^{n} \left( z^{\lambda}_{k} \right)^{d^{\lambda\sigma}_{ki}}.
		\end{align*}
		Hence we obtain $z^{\lambda} \sim z^{\sigma}$.
\end{proof}
\end{prop}

\begin{prop}
\label{prop: hausdorff space}
	The quotient space $X = \bigsqcup_{\lambda \in \Lambda} \mathbb{C}^{n}_{\lambda} / {\sim}$ is a Hausdorff space.
	
\begin{proof}
	Define the projection $\text{pr}: \bigsqcup_{\lambda \in \Lambda} \mathbb{C}^{n}_{\lambda} \to X$ to the quotient space.
	Take two distinct points $[x] \neq [y] \in X$. Let $U_{[x]}$ and $U_{[y]}$ be open subsets containing the points $[x]$ and $[y]$ respectively.
	Then we can write $\text{pr}^{-1}(U_{[x]}),\text{pr}^{-1}(U_{[y]}) \subset X$ as follows:
	\begin{equation*}
		\text{pr}^{-1}(U_{[x]}) = \bigsqcup_{\lambda \in \Lambda}
		U_{[x]}^{\lambda}, \;
		\text{pr}^{-1}(U_{[y]}) = \bigsqcup_{\lambda \in \Lambda}
		U_{[y]}^{\lambda},
	\end{equation*}
	where $U_{[x]}^{\lambda}, U_{[y]}^{\lambda} \subset \mathbb{C}^{n}_{\lambda} \cong \mathbb{C}^{n}$ for $\lambda \in \Lambda$.

	Let $B_{\varepsilon}(x)$ be an open ball of radius $\varepsilon > 0$.
	We define the map $\varphi_{\lambda} : \mathbb{C}_{\lambda}^{n}/ {\sim} \to \mathbb{C}_{\lambda}
	^{n}$ by $\varphi_{\lambda}([z^{\lambda}]) = z^{\lambda}$.
	If $\mathbb{C}^{n}_{\lambda}$ contains the points $x$ and $y$, then there exist $\varepsilon,\varepsilon^{\prime} >0$ such that $B_{\varepsilon}(x) \cap B_{\varepsilon^{\prime}}(y) = \emptyset$. Thus we have $\text{pr}^{-1}(\text{pr}(B_{\varepsilon}(x))) \cap \text{pr}^{-1}(\text{pr}(B_{\varepsilon^{\prime}}(y))) = \emptyset$.
	
	If $x \in \mathbb{C}^{n}_{\lambda}, y \in \mathbb{C}^{n}_{\mu}$ ($\lambda \neq \mu$),
	then there exists an element $\sigma \in \Lambda$ such that 
	$\varphi_{\sigma} \circ \varphi_{\lambda}^{-1}(B_{\varepsilon}(x) \cap U_{\lambda \sigma}) \subset U_{[x]}^{\sigma}$, $\varphi_{\sigma} \circ \varphi_{\mu}^{-1}(B_{\varepsilon^{\prime}}(y) \cap U_{\mu\sigma}) \subset U_{[y]}^{\sigma}$. Thus we can take sufficiently small $\varepsilon,\varepsilon^{\prime} >0$ such that 
	$\text{pr}^{-1}(\text{pr}(\varphi_{\sigma} \circ \varphi_{\lambda}^{-1}(B_{\varepsilon}(x) \cap U_{\sigma\lambda}))) \cap \text{pr}^{-1}(\text{pr}(\varphi_{\sigma} \circ \varphi_{\mu}^{-1}(B_{\varepsilon^{\prime}}(y) \cap U_{\sigma\mu}))) = \emptyset$.

	Suppose that $x \in \mathbb{C}^{n}_{\lambda} \setminus U_{\lambda\mu}, y \in \mathbb{C}^{n}_{\mu} \setminus U_{\mu\lambda}$ ($\lambda \neq \mu$).
	If there exist $\varepsilon,\varepsilon^{\prime} >0$ such that $\mathrm{pr}(B_{\varepsilon}(x)) \cap \mathrm{pr}(B_{\varepsilon^{\prime}}(y)) \neq \emptyset$, then there exist $z_{x} \in B_{\varepsilon}(x) \cap U_{\lambda\mu}$ and $z_{y} \in B_{\varepsilon^{\prime}}(y) \cap U_{\mu\lambda}$ such that $\mathrm{pr}(z_{x}) = \mathrm{pr}(z_{y})$.
	Since $x \not\in U_{\lambda\mu}$, we obtain 
	\begin{equation*}
		0 < \vert x - z_{x} \vert < \varepsilon.
	\end{equation*}
	We can retake $\varepsilon$ smaller than $\vert x - z_{x} \vert$ so that $\mathrm{pr}(B_{\varepsilon}(x)) \cap \mathrm{pr}(B_{\varepsilon^{\prime}}(y)) = \emptyset$.
	Thus we have $\mathrm{pr}^{-1}(\mathrm{pr}(B_{\varepsilon}(x))) \cap \mathrm{pr}^{-1}(\mathrm{pr}(B_{\varepsilon^{\prime}}(y))) = \emptyset$.
	
	Therefore, the quotient space $X$ is a Hausdorff space.
\end{proof}
\end{prop}

Let $U_{\lambda} = \{ [z^{\lambda}] \in X \mid z^{\lambda} \in \mathbb{C}^{n}_{\lambda} \} \subset X$. Then we see $X = \bigcup_{\lambda \in \Lambda} U_{\lambda}$ from Proposition \ref{prop: equiv relation}.
We define a map $\varphi_{\lambda}: U_{\lambda} \to \mathbb{C}^{n}$ by $\varphi_{\lambda} ([z^{\lambda}]) = z^{\lambda}$ for each $\lambda \in \Lambda$.
The following lemma is obvious from the construction above.

\begin{lemma}
\label{lemma: coordinate transformation on inhomogeneous coordinates}
	For all $\lambda, \mu \in \Lambda$ such that $U_{\lambda} \cap U_{\mu} \neq \emptyset$, we have
	\begin{equation}
		\label{eq: coordinate transformation}
		\varphi_{\mu} \circ \varphi_{\lambda}^{-1} (z^{\lambda}) 
		= \left( \prod_{j=1}^{n} \left(z^{\lambda}_{j}\right)^{d^{\lambda\mu}_{j1}}, \ldots, \prod_{j=1}^{n} \left(z^{\lambda}_{j}\right)^{d^{\lambda\mu}_{jn}} \right).
	\end{equation}
\end{lemma}

\begin{defin}
	The set $\{(U_{\lambda}, \varphi_{\lambda})\}_{\lambda \in \Lambda}$ is a system of complex coordinate charts on $X$, whose coordinate transformation is given by Equation \ref{eq: coordinate transformation}.
\end{defin}

\section{Toric Manifolds}
\label{sec: toric manifold}

In this section, we write the inhomogeneous coordinate charts on a toric manifold in terms of the coordinate charts given in Section \ref{sec: coordinates chart}. We also discuss the complements of toric divisors, which we call the \textit{toric divisor complements}.

\subsection{Convex Polytopes and Convex Cones}
\label{subsec: convex polytopes, cones}
We review the definitions and some of the facts of \textit{convex polytopes} and \textit{convex cones} in $\mathbb{R}^{n}$, which are used later.

We first deal with convex polytopes, which are defined as follows:

\begin{defin}
	\label{def: convex polytopes}
	Let $V = \{x_{1},\ldots, x_{s} \} \neq \emptyset$ be a finite set of elements in $\mathbb{R}^{n}$.
	The convex hull $\Delta = \mathrm{conv}(V)$ of $V$ is a convex polytope in $\mathbb{R}^{n}$. Concretely, $\Delta$ is written as 
	\begin{equation*}
		\label{eq: convex polytope def}
		\Delta = \mathrm{conv}(V)
		=
		\left\{ 
			\sum_{i=1}^{s}r_{i}x_{i} \; \middle| \; r_{i} \geq 0, \sum_{i=1}^{s}r_{i}= 1, x_{i} \in V
		\right\}.
	\end{equation*}
	If a convex polytope $\Delta$ can be written as Equation \ref{eq: convex polytope def}, then we say that $\Delta$ is generated by $V = \{x_{1},\ldots, x_{s} \}$.
\end{defin}

The next lemma is obvious.

\begin{lemma}
	\label{lemma: convex polytope linear map}
	Let $f: \mathbb{R}^{n} \to \mathbb{R}^{m}$ be a linear map. 
	If $\Delta \subset \mathbb{R}^{n}$ is a convex polytope generated by $V = \{x_{1},\ldots, x_{s} \}$, then $f(\Delta) \subset \mathbb{R}^{m}$ is also a convex polytope generated by $V_{f} = \{f(x_{1}),\ldots, f(x_{s}) \}$.
\end{lemma}

We deal with convex cones, which are defined as follows:
\begin{defin}
	A subset $\mathcal{C}$ in $\mathbb{R}^{n}$ is a (convex polyhedral) cone if there exist elements $v_{1},\ldots,v_{s} \in \mathcal{C}$ such that 
	\begin{equation}
		\label{eq: def of cones}
		\mathcal{C} = \mathbb{R}_{\geq 0}v_{1} + \cdots + \mathbb{R}_{\geq 0}v_{s}.
	\end{equation}
	If a cone $\mathcal{C}$ can be written as Equation \ref{eq: def of cones}, then we say that $\mathcal{C}$ is generated by $\{v_{1},\ldots,v_{s}\}$. 
\end{defin}

The next lemma is obvious.

\begin{lemma}
	\label{lemma: cone linear map}
	Let $f: \mathbb{R}^{n} \to \mathbb{R}^{m}$ be a linear map. 
	If $\mathcal{C} \subset \mathbb{R}^{n}$ is a cone generated by $\{v_{1},\ldots,v_{s}\}$, then $f(\mathcal{C}) \subset \mathbb{R}^{m}$ is also a cone generated by $\{f(v_{1}),\ldots, f(v_{s})\}$.
\end{lemma}

\begin{defin}
	Let $\mathcal{C} \subset \mathbb{R}^{n}$ be a cone generated by $\{v_{1},\ldots,v_{s}\}$. 
	The point $0 =(0,\ldots,0) \in \mathcal{C} \subset \mathbb{R}^{n}$ is a \textit{vertex} of $\mathcal{C}$ if $\mathcal{C}$ does not contain a nontrivial subspace.
\end{defin}

\begin{lemma}
	Let $\mathcal{C} \subset \mathbb{R}^{n}$ be a cone generated by $\{v_{1},\ldots,v_{s}\}$. 
	The cone $\mathcal{C}$ does not contain a nontrivial subspace 
	if and only if 
	the following is satisfied:
	\begin{equation}
		\label{eq: vertex condition}
		r_{1}v_{1} + \cdots + r_{s} v_{s} = 0, \; r_{i} \geq 0
		\Rightarrow
		r_{1} = \cdots = r_{s} = 0.
	\end{equation}
\begin{proof}
	We first show that if $\mathcal{C}$ does not contain a nontrivial subspace, then Equation \ref{eq: vertex condition} holds.
	We give a proof by showing the contraposition.

	Suppose that $r_{1},\ldots, r_{s} \in \mathbb{R}_{\geq 0}$ satisfy $\sum_{i=1}^{s}r_{i}v_{i} = 0$. 
	Suppose further that there exists some $i_{0} \in \{1,\ldots, s \}$ such that $r_{i_{0}} > 0$.
	Then since we can calculate 
	\begin{align*}
		v_{i{_{0}}} = - \frac{1}{r_{i_{0}}}\sum_{i \neq i_{0}}r_{i}v_{i} 
		= - \sum_{i \neq i_{0}} \frac{r_{i}}{r_{i_{0}}} v_{i},
	\end{align*}
	$W := \{rv_{i_{0}}\mid r \in \mathbb{R} \} \subset \mathcal{C}$ holds.
	Indeed, if $r \geq 0$, then $rv_{i_{0}} \in \mathcal{C}$ by the definition of $\mathcal{C}$; if otherwise, then since from the above calculation we see 
	\begin{equation*}
		rv_{i_{0}} = \sum_{i \neq i_{0}} (-r) \frac{r_{i}}{r_{i_{0}}} v_{i}
	\end{equation*}
	and $(-r) \frac{r_{i}}{r_{i_{0}}} \geq 0$ for any $i \neq i_{0}$, $rv_{i_{0}} \in \mathcal{C}$.
	Since $W$ is a nontrivial subspace in $\mathbb{R}^{n}$, we obtain the contraposition to the desired result.

	We then show that if Equation \ref{eq: vertex condition} holds, then $\mathcal{C}$ does not contain a nontrivial subspace.

	Let $W \neq \emptyset$ be a subspace contained in $\mathcal{C}$.
	Since $W$ is a linear space, if $w\in W$ then $-w \in W$ holds.
	Since $W \subset \mathcal{C}$, there exsit $r_{1},\ldots, r_{s}, r^{\prime}_{1},\ldots, r^{\prime}_{s} \geq 0$ such that 
	\begin{align*}
		w = \sum_{i=1}^{s}r_{i}v_{i},
		\;
		-w = \sum_{i=1}^{s}r^{\prime}_{i}v_{i}.
	\end{align*}
	Since $w + (-w) = 0$, we see that 
	\begin{equation*}
		\sum_{i=1}^{s} (r_{i}+r^{\prime}_{i})v_{i} = 0.
	\end{equation*}
	Since we assume that Equation \ref{eq: vertex condition} holds, $r_{i}+r^{\prime}_{i} = 0$ holds for any $i=1,\ldots,s$.
	Furthermore, since $r_{1},\ldots, r_{s}, r^{\prime}_{1},\ldots, r^{\prime}_{s} \geq 0$, $r_{i} = r^{\prime}_{i} =0$ holds for any $i=1,\ldots,s$.
	This implies that $w = 0$, i.e., $W = \{0\}$.
\end{proof}
\end{lemma}

From the above lemma, we can use the following definition of a \textit{vertex} in a cone.
\begin{defin}
	\label{def: vertex of cone}
	Let $\mathcal{C} \subset \mathbb{R}^{n}$ be a cone generated by $\{v_{1},\ldots,v_{s}\}$. 
	The point $0 =(0,\ldots,0) \in \mathcal{C} \subset \mathbb{R}^{n}$ is a \textit{vertex} of $\mathcal{C}$ if Equation \ref{eq: vertex condition} is satisfied.
\end{defin}

\subsection{An Alternative Construction of Toric Manifolds}
\label{subsec: construction of toric manifolds}
We briefly review the Delzant construction \cite{MR984900} in order to construct inhomogeneous coordinate charts on a toric manifold.
Delzant showed that there is a one-to-one correspondence between compact symplectic toric manifolds and Delzant polytopes, which are moment polytopes for the Hamiltonian torus action on toric manifolds (see \cite[Chapter 1]{MR1301331} for detailed explanations about the Delzant construction).
Delzant polytopes are defined as follows:

\begin{defin}
\label{def: delzant}
	Delzant polytopes are convex polytopes $\Delta$ in $(\mathfrak{t}^{n})^{\ast} \cong \mathbb{R}^{n}$ satisfying the following three conditions:
	\begin{itemize}
		\item simple; each vertex has $n$ edges,
		\item rational; the direction vectors $v^{\lambda}_{1}, \ldots, v^{\lambda}_{n}$ from any vertex $\lambda \in \Lambda$ are integral vectors,
		\item smooth; the vectors $v^{\lambda}_{1}, \ldots, v^{\lambda}_{n}$ chosen as above form a basis of $\mathbb{Z}^{n}$,
	\end{itemize}
	where $\Lambda$ is the set of the vertices in $\Delta$.
\end{defin}
We can define Delzant polytopes in terms of facets in $\Delta$ instead of edges (see \cite[Theorem 4]{MR4417716} for example).
\begin{defin}
\label{def: delzant 2}
	Delzant polytopes are convex polytopes $\Delta$ in $(\mathfrak{t}^{n})^{\ast} \cong \mathbb{R}^{n}$ satisfying the following three conditions:
	\begin{itemize}
			\item simple; each vertex meets $n$ facets,
			\item rational; the inward pointing normal vectors $u^{\lambda}_{1}, \ldots, u^{\lambda}_{n}$ for facets meeting a vertex $\lambda \in \Lambda$ are integral vectors,
			\item smooth; the vectors $u^{\lambda}_{1}, \ldots, u^{\lambda}_{n}$ chosen as above form a basis of $\mathbb{Z}^{n}$,
	\end{itemize}
	where $\Lambda$ is the set of the vertices in $\Delta$.
\end{defin}

We can see that two ways to define Delzant polytopes are equivalent. 
Although we can find a similar result in \cite[Proposition 2.2]{MR2282365}, we give a proof because we shall use the statement repeatedly.
\begin{lemma}
\label{lemma: direction vectors vs normal vectors}
	Let $v^{\lambda}_{1},\ldots,v^{\lambda}_{n}$ be the direction vectors and $u^{\lambda}_{1},\ldots, u^{\lambda}_{n}$ the inward pointing normal vectors for $\lambda \in \Lambda$. 
	Then, 
	\begin{equation*}
		[v^{\lambda}_{1} \cdots v^{\lambda}_{n}] \left[
			\begin{array}{c}
				^t u^{\lambda}_{1} \\
				\vdots \\
				^t u^{\lambda}_{n}
			\end{array}
		\right] = E_{n},
	\end{equation*}
	where $E_{n}$ denotes the identity matrix.

\begin{proof}
	Let $e_{1}, \ldots, e_{n}$ be the standard basis in $\mathbb{R}^{n}$.
	Since the direction vectors $v^{\lambda}_{1}, \ldots, v^{\lambda}_{n}$ form a basis of $\mathbb{Z}^{n}$, there exists a square matrix $B_{n}^{\lambda}$ such that 
	\begin{equation}
		\label{eq: direction vectors}
		E_{n} = [e_{1} \cdots e_{n} ]
		= [ v^{\lambda}_{1} \cdots v^{\lambda}_{n} ] B_{n}^{\lambda}.
	\end{equation}
	Since the matrix $B_{n}^{\lambda}$ is the inverse matrix for the matrix $[ v^{\lambda}_{1} \cdots v^{\lambda}_{n} ]$, we see $B_{n}^{\lambda}[ v^{\lambda}_{1} \cdots v^{\lambda}_{n} ] = E_{n}$.
	Moreover, the matrix $B^{\lambda}_{n}$ is in $GL(n;\mathbb{Z})$ from (\ref{eq: direction vectors}).

	We define the vectors $u^{\lambda}_{1}, \ldots, u^{\lambda}_{n} \in \mathbb{Z}^{n}$ by 
	\begin{equation*}
		B_{n}^{\lambda} = \left[
			\begin{array}{c}
				^t u^{\lambda}_{1} \\
				\vdots \\
				^t u^{\lambda}_{n}
			\end{array}
		\right].
	\end{equation*}
	By calculating $B_{n}^{\lambda} [ v^{\lambda}_{1} \cdots v^{\lambda}_{n} ]$, we obtain 
	\begin{align*}
		B_{n}^{\lambda} [ v^{\lambda}_{1} \cdots v^{\lambda}_{n} ]
		&= \left[
			\begin{array}{c}
				^t u^{\lambda}_{1} \\
				\vdots \\
				^t u^{\lambda}_{n}
			\end{array}
		\right] [ v^{\lambda}_{1} \cdots v^{\lambda}_{n} ] \\
		&= \left[
			\begin{array}{ccc}
				^t u^{\lambda}_{1}v^{\lambda}_{1} & \cdots & ^t u^{\lambda}_{1}v^{\lambda}_{n} \\
				\vdots & \ddots & \vdots \\
				^t u^{\lambda}_{n}v^{\lambda}_{1} & \cdots & ^t u^{\lambda}_{n}v^{\lambda}_{n}
			\end{array}
		\right] \\
		&= \left[
			\begin{array}{ccc}
				\langle u^{\lambda}_{1}, v^{\lambda}_{1} \rangle & \cdots & \langle u^{\lambda}_{1}, v^{\lambda}_{n} \rangle \\
				\vdots & \ddots & \vdots \\
				\langle u^{\lambda}_{n}, v^{\lambda}_{1} \rangle & \cdots & \langle u^{\lambda}_{1}, v^{\lambda}_{n} \rangle 
			\end{array}
		\right].
	\end{align*}
	Since $B_{n}^{\lambda}[ v^{\lambda}_{1} \cdots v^{\lambda}_{n} ] = E_{n}$, we have $\langle u^{\lambda}_{i}, v^{\lambda}_{j} \rangle = \delta_{ij}$ (Kronecker's delta). We say that 
	\begin{equation*}
		u^{\lambda}_{1} \in \mathrm{span}\{v^{\lambda}_{2},\ldots, v^{\lambda}_{n}\}^{\perp},
		\ldots,
		u^{\lambda}_{n} \in \mathrm{span}\{v^{\lambda}_{1},\ldots, v^{\lambda}_{n-1}\}^{\perp},
	\end{equation*}
	which means that the vectors $u^{\lambda}_{1}, \ldots, u^{\lambda}_{n}$ are inward pointing normal vectors to facets meeting the vertex $\lambda$.
\end{proof}
\end{lemma}

For $\lambda \in \Lambda$, we define an $n \times n$ matrix $Q^{\lambda} = \left[ v^{\lambda}_{1}, \ldots, v^{\lambda}_{n} \right] (= \left[ Q^{\lambda}_{ij} \right])$. In general $\det Q^{\lambda} = \pm 1$ by the definition, but we assume $\det Q^{\lambda} = 1$ by changing the numbering of $v^{\lambda}_{1}, \ldots, v^{\lambda}_{n}$.
We also define a matrix $D^{\lambda \mu}$ by $D^{\lambda \mu} = (Q^{\lambda})^{-1} Q^{\mu} (= [d^{\lambda\mu}_{ij}])$ for each $\lambda, \mu \in \Lambda$ as we defined in Section \ref{sec: coordinates chart}.

From the construction in Section \ref{sec: coordinates chart}, we obtain a system of complex coordinate charts $\{(U_{\lambda}, \varphi_{\lambda})\}_{\lambda \in \Lambda}$ on a toric manifold $X$ associated with a Delzant polytope $\Delta$.

\begin{remark}
Azam, Cannizzo, and Lee explained the construction of symplectic toric manifolds with a system of the inhomogeneous coordinate charts from a data of Delzant polytopes \cite{MR4417716}.
In this case, 
the coordinate transformation of 
our system of complex coordinate charts $\{(U_{\lambda}, \varphi_{\lambda})\}_{\lambda \in \Lambda}$
coincides with the one constructed in \cite{MR4417716}.
Hereafter, we call $\{(U_{\lambda}, \varphi_{\lambda})\}_{\lambda \in \Lambda}$ a system of the inhomogeneous coordinate chart on a toric manifold.
\end{remark}

From this section, 
we write $X$ by a compact toric manifold of complex dimension $n$, 
$\Delta$ by the Delzant polytope of $X$, 
$\Lambda$ by the set of the vertices in the polytope $\Delta$.

\begin{remark}
The coordinate transformation of the inhomogeneous coordinates also coincides with the one in algebraic geometry (see for example \cite{MR1234037}).
Note that we may have a fan of toric manifolds $X$ by taking integral vectors inward pointing normal to each facets of Delzant polytopes of $X$.
\end{remark}

\subsection{Toric Divisor Complements}
\label{subsec: toric divisor complements}

In this section, we construct a diffeomorphism between the complements of toric divisors in $X$ and $(\mathbb{C}^{\ast})^{n}$.

We define $\check{U}_{\lambda} = \{[z^{\lambda}] \in U_{\lambda} \mid z^{\lambda}_{1}z^{\lambda}_{2} \cdots z^{\lambda}_{n} \neq 0\} \subset X$ for each $\lambda \in \Lambda$, then we have $\check{M} = \bigcup_{\lambda \in \Lambda} \check{U}_{\lambda}$. 
We call $\check{M}$ be the \textit{toric divisor complement}. 
Furthermore, we see $\check{U}_{\lambda} = \check{U}_{\sigma} = \check{M}$ for $\lambda, \sigma \in \Lambda$.

\begin{defin}
	We define a map $\phi_{\lambda}: \varphi_{\lambda} (\check{U}_{\lambda}) \to (\mathbb{C}^{\ast})^{n}$ by 
	\begin{equation*}
		\phi_{\lambda} (z^{\lambda}_{1}, \ldots, z^{\lambda}_{n}) = \left(\prod_{j=1}^{n} (z^{\lambda}_{j})^{\hat{Q}^{\lambda}_{j1}}, \ldots, \prod_{j=1}^{n} (z^{\lambda}_{j})^{\hat{Q}^{\lambda}_{jn}} \right),
	\end{equation*}
	where $(Q^{\lambda})^{-1} = [\hat{Q}^{\lambda}_{ij}]$.
\end{defin}

\begin{lemma}
\label{lemma: toric divisor complement to complex torus}
	For any $\lambda, \sigma \in \Lambda$, 
	$\phi_{\lambda}\circ \varphi_{\lambda} = \phi_{\sigma} \circ \varphi_{\sigma}$.

\begin{proof}
	Since we define $D^{\sigma\lambda} = (Q^{\sigma})^{-1} Q^{\lambda}$, we see $D^{\sigma\lambda}(Q^{\lambda})^{-1} = (Q^{\sigma})^{-1}$.
	For any $[z^{\lambda}] \in \check{U}_{\lambda}$, we have 
	\begin{align*}
		\phi_{\lambda}\circ \varphi_{\lambda} ([z^{\lambda}])
		=& \phi_{\lambda}\circ \varphi_{\lambda}\left( \left[\left( \prod_{j=1}^{n} (z^{\sigma}_{j})^{d^{\sigma\lambda}_{j1}}, \ldots, \prod_{j=1}^{n} (z^{\sigma}_{j})^{d^{\sigma\lambda}_{jn}}\right) \right] \right) \\
		=& \phi_{\lambda} \left( \prod_{j=1}^{n} (z^{\sigma}_{j})^{d^{\sigma\lambda}_{j1}}, \ldots, \prod_{j=1}^{n} (z^{\sigma}_{j})^{d^{\sigma\lambda}_{jn}} \right) \\
		=& \left( \prod_{i=1}^{n} \left(\prod_{j=1}^{n} (z^{\sigma}_{j})^{d^{\sigma\lambda}_{ji}} \right)^{\hat{Q}^{\lambda}_{i1}}, \ldots, \prod_{i=1}^{n} \left( \prod_{j=1}^{n} (z^{\sigma}_{j})^{d^{\sigma\lambda}_{ji}} \right)^{\hat{Q}^{\lambda}_{in}} \right) \\
		=& \left( \prod_{j=1}^{n} (z^{\sigma}_{j})^{\hat{Q}^{\sigma}_{j1}}, \ldots, \prod_{j=1}^{n} (z^{\sigma}_{j})^{\hat{Q}^{\sigma}_{jn}} \right) \\
		=& \phi_{\sigma} \circ \varphi_{\sigma} ([z^{\sigma}]).
	\end{align*}
\end{proof}
\end{lemma}
From Lemma \ref{lemma: toric divisor complement to complex torus}, we can define the following map independent of the choice of $\lambda \in \Lambda$.

\begin{defin}
\label{def: toric divisor complement to complex torus}
We define a map $\phi: \check{M} \to (\mathbb{C}^{\ast})^{n} = \{(z_{1}, \ldots, z_{n})\mid z_{1} z_{2} \cdots z_{n} \neq 0\}$ by $\phi = \phi_{\lambda} \circ \varphi_{\lambda}$.
\end{defin}

Next we construct the inverse map $\hat{\phi}:(\mathbb{C}^{\ast})^{n} \to \check{M}$, which is actually similar to the construction of $\phi$.

\begin{defin}
	\label{def: inverse map of toric divisor complement to complex torus}
	We define a map $\hat{\phi}_{\lambda}: (\mathbb{C}^{\ast})^{n} \to \varphi_{\lambda} (\check{U}_{\lambda})$ by 
	\begin{equation*}
		\hat{\phi}_{\lambda} (z_{1}, \ldots, z_{n}) = \left(\prod_{j=1}^{n}(z_{j})^{Q^{\lambda}_{j1}}, \ldots, \prod_{j=1}^{n} (z_{j})^{Q^{\lambda}_{jn}} \right).
	\end{equation*}
\end{defin}

\begin{lemma}
	\label{lemma: complex torus to toric divisor complement}
	For any $\lambda, \sigma \in \Lambda$, $\varphi_{\lambda}^{-1} \circ \hat{\phi}_{\lambda} = \varphi_{\sigma}^{-1} \circ \hat{\phi}_{\sigma}$.

\begin{proof}
	Since we define $D^{\lambda\sigma} = (Q^{\lambda})^{-1} Q^{\sigma}$, we see $Q^{\lambda}D^{\lambda\sigma} = Q^{\sigma}$.
	For any $z = (z_{1}, \ldots, z_{n}) \in (\mathbb{C}^{\ast})^{n}$, we have 
	\begin{align*}
		\varphi_{\lambda}^{-1} \circ \hat{\phi}_{\lambda} (z) 
		=& \varphi_{\lambda}^{-1}\left(\prod_{j=1}^{n}(z_{j})^{Q^{\lambda}_{j1}}, \ldots, \prod_{j=1}^{n} (z_{j})^{Q^{\lambda}_{jn}} \right) \\
		=& \left[\left(\prod_{j=1}^{n}(z_{j})^{Q^{\lambda}_{j1}}, \ldots, \prod_{j=1}^{n} (z_{j})^{Q^{\lambda}_{jn}} \right)\right] \\
		=& \left[ \left( \prod_{i=1}^{n} \left(\prod_{j=1}^{n}(z_{j})^{Q^{\lambda}_{ji}} \right)^{d^{\lambda\sigma}_{i1}}, \ldots, \prod_{i=1}^{n} \left( \prod_{j=1}^{n} (z_{j})^{Q^{\lambda}_{ji}} \right)^{d^{\lambda\sigma}_{in}} \right)\right] \\
		=& \left[\left( \prod_{j=1}^{n} (z_{j})^{Q^{\sigma}_{j1}}, \ldots, \prod_{j=1}^{n} (z_{j})^{Q^{\sigma}_{jn}} \right)\right] \\
		=&\varphi_{\sigma}^{-1} \circ \hat{\phi}_{\sigma} (z).
	\end{align*}
\end{proof}
\end{lemma}

From Lemma \ref{lemma: complex torus to toric divisor complement}, we can define the following map independent of the choice of $\lambda \in \Lambda$.

\begin{defin}
\label{def: complex torus to toric divisor complement}
	We define a map $\hat{\phi}:(\mathbb{C}^{\ast})^{n} \to \check{M}$ by $\hat{\phi} = \varphi_{\lambda}^{-1} \circ \hat{\phi}_{\lambda}$.
\end{defin}

We check that the map $\hat{\phi}$ defined in Definition \ref{def: complex torus to toric divisor complement} is the inverse map of $\phi$ defined in Definition \ref{def: toric divisor complement to complex torus}.

\begin{lemma}
\label{lemma: inverse map}
	$\hat{\phi}\circ \phi = \text{id}_{\check{M}}$, $\phi \circ \hat{\phi} = \text{id}_{(\mathbb{C}^{\ast})^{n}}$.

\begin{proof}
	Since we define $\phi = \phi_{\lambda} \circ \varphi_{\lambda}$ and $\hat{\phi} = \varphi_{\lambda}^{-1} \circ \hat{\phi}_{\lambda}$, we obtain 
	$\hat{\phi}\circ \phi = \varphi_{\lambda}^{-1}\circ \hat{\phi}_{\lambda} \circ \phi_{\lambda} \circ \varphi_{\lambda}$ and $\phi \circ \hat{\phi} = \phi_{\lambda} \circ \varphi_{\lambda} \circ \varphi_{\lambda}^{-1} \circ \hat{\phi}_{\lambda} = \phi_{\lambda} \circ \hat{\phi}_{\lambda}$.
	We say that it is sufficient to show that $\hat{\phi}_{\lambda} \circ \phi_{\lambda} = \text{id}_{\varphi_{\lambda}(\check{U}_{\lambda})}$ and $\phi_{\lambda} \circ \hat{\phi}_{\lambda} = \text{id}_{(\mathbb{C}^{\ast})^{n}}$.

	For $z^{\lambda} \in \varphi_{\lambda}(\check{U}_{\lambda})$, by similar calculation, we have 
	\begin{align*}
		\hat{\phi}_{\lambda} \circ \phi_{\lambda} (z^{\lambda}_{1}, \ldots, z^{\lambda}_{n})
		=& \hat{\phi}_{\lambda} \left(\prod_{j=1}^{n}(z^{\lambda}_{j})^{\hat{Q}^{\lambda}_{j1}}, \ldots, \prod_{j=1}^{n}(z^{\lambda}_{j})^{\hat{Q}^{\lambda}_{jn}} \right) \\
		=& \left( \prod_{i=1}^{n} \left(\prod_{j=1}^{n}(z^{\lambda}_{j})^{\hat{Q}^{\lambda}_{ji}} \right)^{Q^{\lambda}_{i1}}, \ldots, \prod_{i=1}^{n}\left(\prod_{j=1}^{n}(z^{\lambda}_{j})^{\hat{Q}^{\lambda}_{ji}} \right)^{Q^{\lambda}_{in}} \right) \\
		=& \left( \prod_{j=1}^{n}(z^{\lambda}_{j})^{\delta_{j1}}, \ldots, \prod_{j=1}^{n}(z^{\lambda}_{j})^{\delta_{jn}} \right) \\
		=& (z^{\lambda}_{1}, \ldots, z^{\lambda}_{n}).
	\end{align*}
	Thus we obtain $\hat{\phi}_{\lambda} \circ \phi_{\lambda} = \text{id}_{\varphi_{\lambda}(\check{U}_{\lambda})}$. 

	For $z \in (\mathbb{C}^{\ast})^{n}$, by similar calculation, we have 
	\begin{align*}
		\phi_{\lambda} \circ \hat{\phi}_{\lambda} (z_{1},\ldots, z_{n}) 
		=& \phi_{\lambda} \left( \prod_{j=1}^{n}(z_{j})^{Q^{\lambda}_{j1}}, \ldots, \prod_{j=1}^{n}(z_{j})^{Q^{\lambda}_{jn}} \right) \\
		=& \left( \prod_{i=1}^{n}\left(\prod_{j=1}^{n}(z_{j})^{Q^{\lambda}_{ji}} \right)^{\hat{Q}^{\lambda}_{i1}}, \ldots, \prod_{i=1}^{n}\left(\prod_{j=1}^{n}(z_{j})^{Q^{\lambda}_{ji}} \right)^{\hat{Q}^{\lambda}_{in}} \right) \\
		=& (z_{1}, \ldots, z_{n}).
	\end{align*}
	Thus we obtain $\phi_{\lambda} \circ \hat{\phi}_{\lambda} = \text{id}_{(\mathbb{C}^{\ast})^{n}}$.
\end{proof}
\end{lemma}

\section{Torus-equivariantly Embedded Toric Manifolds}
\label{sec: toric submanifolds}

We construct $k$-dimensional complex submanifolds $\overline{C(V)}$ in toric manifolds $X$ associated to affine subspaces $V$ in $\mathbb{R}^{n} \cong \mathfrak{t}^{n}$ and examine their fundamental properties.

In Section \ref{subsec: construction of toric submanifolds}, we give the construction of $\overline{C(V)}$.
In Section \ref{subsec: subtorus actions on toric manifolds}, we consider a Hamiltonian subtorus action on a toric manifold $X$.
In Section \ref{subsec: torus actions on toric submanifolds}, we consider the Hamiltonian torus action on $\overline{C(V)}$.

\subsection{Construction of Torus-equivariantly Embedded Toric Manifolds}
\label{subsec: construction of toric submanifolds}

First, we will concentrate on Yamamoto's construction \cite[Lemma 6.1]{MR3856842} of complex submanifolds $C(V)$ in $(\mathbb{C}^{\ast})^{n}$.
Let $e_{1}, \ldots, e_{n}$ be the standard basis in $\mathbb{R}^{n}$.
We write $\langle, \rangle$ for the inner product of vectors.
Fix $k = 1, \ldots, n$.
Let $p_{1}, \ldots, p_{k} \in \mathbb{Z}^{n}$ be primitive vectors which are linearly independent, and $a \in \mathbb{R}^{n}$.
Then we consider an affine subspace $V = \mathbb{R}p_{1} + \cdots + \mathbb{R}p_{k} + a \cong \mathbb{R}^{k}$ in $\mathbb{R}^{n}$, which may have rational slope.
Yamamoto constructed $k$-dimensional complex manifolds $C(V)$ in $(\mathbb{C}^{\ast})^{n}$. 
Although we do not give the same statement as the original one, the statement is like as follows.

\begin{prop}
\label{prop: complex submanifolds in complex torus}
Given an affine subspace $V = \mathbb{R}p_{1} + \cdots + \mathbb{R}p_{k} + a$ in $\mathbb{R}^{n}$, we can construct a complex submanifold $C(V) \cong (\mathbb{C}^{\ast})^{k}$ in $\check{M} \cong (\mathbb{C}^{\ast})^{n}$ by 
\begin{equation*}
	C(V) = \left\{ (e^{x_{1}+\sqrt{-1}y_{1}}, \ldots, e^{x_{n}+ \sqrt{-1}y_{n}} ) \in (\mathbb{C}^{\ast})^{n} \;\middle|\; 
	\begin{aligned} 
		x_{i} &= \sum_{l=1}^{k} \langle p_{l}, e_{i} \rangle u_{l} + \langle a, e_{i} \rangle, \\
		y_{i} &= \sum_{l=1}^{k} \langle p_{l}, e_{i} \rangle v_{l}
	\end{aligned}
	\right\},
\end{equation*}
where $(e^{u_{1}+\sqrt{-1}v_{1}}, \ldots, e^{u_{k}+\sqrt{-1}v_{k}}) \in (\mathbb{C}^{\ast})^{k}$.
\end{prop}

Note that if $k=0$, then $C(V)$ is a point $(e^{\langle a,e_{1} \rangle}, \ldots, e^{\langle a,e_{n} \rangle})$ in $(\mathbb{C}^{\ast})^{n}$.
We rewrite the expression of $C(V)$ in 
Proposition \ref{prop: complex submanifolds in complex torus} as follows:

\begin{prop}
\label{prop: complex submanifolds in complex torus2}
Let $C(V)$ be a complex submanifold in $\check{M}$ given in Proposition \ref{prop: complex submanifolds in complex torus}. 
There exists a primitive basis $q_{k+1}, \ldots, q_{n} \in \mathbb{Z}^{n}$ of the orthogonal subspace to $V = \mathbb{R}p_{1} + \cdots + \mathbb{R}p_{k} + a$ in $\mathbb{R}^{n}$ such that 
\begin{equation*}
	C(V) = \left\{ (e^{w_{1}}, \ldots, e^{w_{n}}) \in (\mathbb{C}^{\ast})^{n} \;\middle|\; ~^t[q_{k+1} \cdots q_{n}] \left( \left[ \begin{array}{c} w_{1} \\ \vdots \\ w_{n} \end{array} \right] -a \right) = 0 \right\},
\end{equation*}
where $w_{1} = x_{1} + \sqrt{-1}y_{1}, \ldots, w_{n} = x_{n} + \sqrt{-1}y_{n}$.

\begin{proof}
Hereafter, we calculate angle coordinates $y = (y_{1},\ldots, y_{n})$ and $v = (v_{1}, \ldots, v_{k})$ up to $2\pi \mathbb{Z}$. 

From Proposition \ref{prop: complex submanifolds in complex torus}, we have 
\begin{equation*}
	\left[ \begin{array}{c} x_{1} \\ \vdots \\ x_{n} \end{array} \right]
	= \left[p_{1} \cdots p_{k} \right] \left[ \begin{array}{c} u_{1} \\ \vdots \\ u_{k} \end{array} \right] + a, \;
	\left[ \begin{array}{c} y_{1} \\ \vdots \\ y_{n} \end{array} \right]
	= \left[p_{1} \cdots p_{k} \right] \left[ \begin{array}{c} v_{1} \\ \vdots \\ v_{k} \end{array} \right].
\end{equation*}
We obtain 
\begin{equation*}
	\left[ \begin{array}{c} w_{1} \\ \vdots \\ w_{n} \end{array} \right] - a
	= \left[ \begin{array}{c} x_{1}+\sqrt{-1}y_{1} \\ \vdots \\ x_{n}+\sqrt{-1}y_{n} \end{array} \right] -a
	= \left[p_{1} \cdots p_{k} \right] \left[ \begin{array}{c} u_{1}+\sqrt{-1}v_{1} \\ \vdots \\ u_{k}+\sqrt{-1}v_{k} \end{array} \right].
\end{equation*}
We can take a basis of primitive vectors $q_{k+1},\ldots,q_{n} \in \mathbb{Z}^{n}$ of the orthogonal subspace to $V$.
Multiplying both sides by $^t[q_{k+1} \cdots q_{n}]$, we obtain 
\begin{align*}
	& ~^t[q_{k+1} \cdots q_{n}] \left( \left[ \begin{array}{c} w_{1} \\ \vdots \\ w_{n} \end{array} \right] -a \right) \\
	=& \left[ \begin{array}{c} ^tq_{k+1} \\ \vdots \\ ^tq_{n} \end{array} \right] \left[p_{1} \cdots p_{k} \right] \left[ \begin{array}{c} u_{1}+\sqrt{-1}v_{1} \\ \vdots \\ u_{k}+\sqrt{-1}v_{k} \end{array} \right] \\
	=& \left[ 
		\begin{array}{cccc}
			\langle q_{k+1}, p_{1} \rangle & \langle q_{k+1}, p_{2} \rangle & \cdots & \langle q_{k+1}, p_{k} \rangle \\
			\langle q_{k+2}, p_{1} \rangle & \langle q_{k+2}, p_{2} \rangle & \cdots & \langle q_{k+2}, p_{k} \rangle \\
			\vdots & \vdots & \ddots & \vdots \\
			\langle q_{n}, p_{1} \rangle & \langle q_{n}, p_{2} \rangle & \cdots & \langle q_{n}, p_{k} \rangle
		\end{array}
	\right]
	\left[ \begin{array}{c} u_{1}+\sqrt{-1}v_{1} \\ \vdots \\ u_{k}+\sqrt{-1}v_{k} \end{array} \right] \\
	=& 0.
\end{align*}
\end{proof}
\end{prop}

Note that this submanifold $C(V)$ can be regarded as a complex subtorus $(\mathbb{C}^{\ast})^{k}$ in $(\mathbb{C}^{\ast})^{n}$.

We will write $C(V)$ explicitly as a submanifold in $\check{M}$.
Recall that $\check{M} = \bigcup_{\lambda \in \Lambda} \check{U}_{\lambda}$. Using the map $\hat{\phi}_{\lambda} :(\mathbb{C}^{\ast})^{n} \to \varphi_{\lambda}(\check{U}_{\lambda})$, we obtain 
\begin{equation*}
	\hat{\phi}_{\lambda}(C(V)) = 
	\{ z^{\lambda} = e^{w^{\lambda}} \in \varphi_{\lambda}(\check{U}_{\lambda}) \mid ~^t[q_{k+1} \cdots q_{n}] ((Q^{\lambda})^{-1}w^{\lambda} - a) = 0 \},
\end{equation*}
where $z^{\lambda} = e^{w^{\lambda}}$ means that $(z^{\lambda}_{1}, \ldots, z^{\lambda}_{n}) = (e^{w^{\lambda}_{1}}, \ldots, e^{w^{\lambda}_{n}})$.

Next we take the closure of $\hat{\phi}_{\lambda}(C(V)) (\subset \varphi_{\lambda}(U_{\lambda}) \cong \mathbb{C}^{n})$. 
Since $\hat{Q}^{\lambda}_{il} = \langle u^{\lambda}_{i}, e_{l} \rangle$, we have 
\begin{equation*}
	\sum_{l=1}^{n} \hat{Q}^{\lambda}_{il} \langle q_{j},e_{l} \rangle 
	= \langle u^{\lambda}_{i}, q_{j} \rangle
\end{equation*}
for $i=1,\ldots, n$ and $j = k+1, \ldots, n$.
Define three subsets $\mathcal{I}^{+}_{\lambda,j},\mathcal{I}^{-}_{\lambda,j},\mathcal{I}^{0}_{\lambda,j} \subset \{1,2,\ldots, n \}$ by 
\begin{align}
	\label{eq: plus}
	\mathcal{I}^{+}_{\lambda,j} =& \left\{ i \in \{1,2,\ldots, n \} \;\middle|\; \langle u^{\lambda}_{i}, q_{j} \rangle \geq 0 \right\}, \\ 
	\label{eq: minus}
	\mathcal{I}^{-}_{\lambda,j} =& \left\{ i \in \{1,2,\ldots, n \} \;\middle|\; \langle u^{\lambda}_{i}, q_{j} \rangle \leq 0 \right\}, \\
	\label{eq: zero}
	\mathcal{I}^{0}_{\lambda,j} =& \left\{ i \in \{1,2,\ldots, n \} \;\middle|\; \langle u^{\lambda}_{i}, q_{j} \rangle = 0 \right\},
\end{align}
for $\lambda \in \Lambda$ and $j = k+1
,\ldots, n$. 
Note that $\mathcal{I}^{+}_{\lambda,j} \cap \mathcal{I}^{-}_{\lambda,j} = \mathcal{I}^{0}_{\lambda,j}$ and $\mathcal{I}^{+}_{\lambda,j} \cup \mathcal{I}^{-}_{\lambda,j} = \{1,2,\ldots, n\}$ for any $\lambda \in \Lambda$ and $j = k+1, \ldots, n$. 

From the expression of $\hat{\phi}_{\lambda}(C(V))$, direct calculation gives us 
\begin{align*}
	\langle q_{l}, (Q^{\lambda})^{-1}w^{\lambda} - a \rangle
	&=
	\langle q_{l}, (Q^{\lambda})^{-1}\log z^{\lambda} - a \rangle \\
	&=
	\langle {^t (Q^{\lambda})^{-1}} q_{l}, \log z^{\lambda} \rangle - \langle q_{l}, a \rangle \\
	&=
	\left\langle 
	\left[
		\begin{matrix}
			\langle u^{\lambda}_{1}, q_{l} \rangle \\
			\vdots \\
			\langle u^{\lambda}_{n}, q_{l} \rangle
		\end{matrix}
	\right], \log z^{\lambda} \right\rangle 
	- \langle a, q_{l} \rangle \\
	&=
	\log \left( \prod_{j=1}^{n} (z^{\lambda}_{j})^{\langle u^{\lambda}_{j},q_{l} \rangle} \right)
	- \langle a,q_{l} \rangle
\end{align*}
for $l = k+1,\ldots,n$.
We can define $\overline{C(V)} = \bigcup_{\lambda \in \Lambda}\varphi_{\lambda}^{-1}(\overline{C_{\lambda}(V)}) \subset X$ by 
\begin{equation}
	\label{eq: toric submanifolds}
	\overline{C_{\lambda}(V)} = \left\{ z^{\lambda} \in \varphi_{\lambda} (U_{\lambda}) \;\middle|\; f^{\lambda}_{j}(z^{\lambda}) = 0 ~\text{for $j = k+1, \ldots, n$} \right\},
\end{equation}
where $f^{\lambda}_{j}$ is defined by 
\begin{equation}
	\label{eq: defining eqautions for toric submanifolds}
	f^{\lambda}_{j}(z^{\lambda}) 
	= \prod_{i \in \mathcal{I}^{+}_{\lambda,j}} (z^{\lambda}_{i})^{\langle u^{\lambda}_{i}, q_{j} \rangle}
	- e^{\langle a, q_{j} \rangle}\prod_{i \in \mathcal{I}^{-}_{\lambda,j}} (z^{\lambda}_{i})^{-\langle u^{\lambda}_{i}, q_{j} \rangle} 
\end{equation}
for each $j = k+1, \ldots, n$.
Here, if $\mathcal{I}^{+}_{\lambda,j} = \emptyset$, then $\prod_{i \in \mathcal{I}^{+}_{\lambda,j}}(z^{\lambda}_{i})^{\langle u^{\lambda}_{i}, q_{j} \rangle} = 1$. Similarly, if $\mathcal{I}^{-}_{\lambda,j} = \emptyset$, then $\prod_{i \in \mathcal{I}^{-}_{\lambda,j}}(z^{\lambda}_{i})^{-\langle u^{\lambda}_{i}, q_{j} \rangle} = 1$.

$\overline{C_{\lambda}(V)}$ is a zero locus of $f^{\lambda}_{k+1},\ldots, f^{\lambda}_{n}$.
Note that if $k = 0$, then $\overline{C(V)}$ is a point in $X$.

\begin{remark}
\label{remark: toric submanifold}
By the implicit function theorem, $\overline{C(V)}$ is a complex submanifold in $X$ if the rank of the Jacobian matrix of $f^{\lambda}_{k+1}, \ldots, f^{\lambda}_{n}$ is equal to $n-k$ for any points $p \in \overline{C(V)}$ and any $\lambda \in \Lambda$.
\end{remark}

We demonstrate some examples for complex submanifolds $\overline{C(V)}$. 
Example \ref{eg: cp2, o.k.} gives an example for a complex submanifold in $X = \mathbb{C}P^{2}$, 
while Example \ref{eg: cp2, n.g.} deals with a subset in $X = \mathbb{C}P^{2}$ which does not become a complex submanifold in $\mathbb{C}P^{2}$.
In the following two examples,
define the points $\lambda$, $\mu$, $\sigma$ in $(\mathfrak{t}^{2})^{\ast} \cong \mathbb{R}^{2}$ by 
	\begin{equation*}
		\lambda = (0,0), ~
		\mu = (2,0), ~ 
		\sigma = (0,2).
	\end{equation*}
Let $\Delta$ be a polytope defined by the convex hull of the points $\lambda$, $\mu$, $\sigma$.
From the Delzant polytope $\Delta$ of $\mathbb{C}P^{2}$  we define 
the inward pointing normal vectors to the facets by 
\begin{equation*}
	u^{\lambda}_{1} =
	\left[
		\begin{matrix}
			1 \\ 0
		\end{matrix}
	\right],
	u^{\lambda}_{2} =
	\left[
		\begin{matrix}
			0 \\ 1 
		\end{matrix}
	\right], 
	u^{\mu}_{1} = 
	\left[
		\begin{matrix}
			-1 \\ -1 
		\end{matrix}
	\right],
	u^{\mu}_{2} = 
	\left[
		\begin{matrix}
			1 \\ 0
		\end{matrix}
	\right],
	u^{\sigma}_{1} = 
	\left[
		\begin{matrix}
			0 \\ 1 
		\end{matrix}
	\right],
	u^{\sigma}_{2} = 
	\left[
		\begin{matrix}
			-1 \\ -1
		\end{matrix}
	\right],
\end{equation*}
where $\Lambda = \{ \lambda,\mu, \sigma\}$.

\begin{example}
\label{eg: cp2, o.k.}
	Let $X = \mathbb{C}P^{2}$, $k = 1$, $V$ be spanned by $p = {^t[1 ~ 1]}$.
	We can choose a basis $q$ of the orthogonal subspace to $V$ as $q = {^t[1 ~ -1]}$. 
	In this case, the subset $\overline{C(V)}$ is a complex submanifold in $X$.
	Indeed, we give $f^{\lambda},f^{\mu}, f^{\sigma}$ by
	\begin{equation*}
		f^{\lambda} = z^{\lambda}_{1} - z^{\lambda}_{2}, ~
		f^{\mu} = z^{\mu}_{2} - 1, ~
		f^{\sigma} = 1- z^{\sigma}_{1},
	\end{equation*}
	respectively. 
	Since the Jacobian matrices are expressed as 
	\begin{equation*}
		Df^{\lambda} = [1 ~ -1],~
		Df^{\mu} = [0 ~ 1],~
		Df^{\sigma} = [1 ~ 0],
	\end{equation*}
	respectively, we see the rank of each matrix is one.
\end{example}

\begin{example}
\label{eg: cp2, n.g.}
	Let $X = \mathbb{C}P^{2}$, $k = 1$, $V$ be spanned by $p = {^t[3 ~ 1]}$.
	We can choose a basis $q$ of the orthogonal subspace to $V$ as $q = {^t[1 ~ -3]}$. 
	In this case, the subset $\overline{C(V)}$ is not a complex submanifold in $X$.
	Indeed, we give $f^{\lambda},f^{\mu}, f^{\sigma}$ by
	\begin{equation*}
		f^{\lambda} = z^{\lambda}_{1} - (z^{\lambda}_{2})^{3}, ~
		f^{\mu} = (z^{\mu}_{1})^{2}z^{\mu}_{2} - 1, ~
		f^{\sigma} = (z^{\sigma}_{1})^{3} - (z^{\sigma}_{2})^{2},
	\end{equation*}
	respectively. Note that $(0,0) \notin \overline{C_{\mu}(V)}$. 
	Since the Jacobian matrices are expressed as 
	\begin{equation*}
		Df^{\lambda} =  [1 ~ -3(z^{\lambda}_{2})^{2}],~
		Df^{\mu} = [2z^{\mu}_{1}z^{\mu}_{2} ~ (z^{\mu}_{1})^{2}],~
		Df^{\sigma} = [3(z^{\sigma}_{1})^{2} ~ -2z^{\sigma}_{2}],
	\end{equation*}
	respectively, we see the rank of $Df^{\lambda}$ and $Df^{\mu}$ is one.
	However, the rank of $Df^{\sigma}$ becomes zero at the point $(z^{\sigma}_{1},z^{\sigma}_{2})=(0,0) \in \overline{C_{\sigma}(V)}$.
\end{example}

Notice that if we fix a toric manifold $X$, then we may classify examples of complex submanifolds $\overline{C(V)}$ in $X$ in terms of the conditions of $V$.
Other examples for $X = \mathbb{C}P^{2}$ are treated in Section \ref{subsec: example}.

Suppose that $\overline{C(V)}$ is a complex submanifold in $X$.
Then, there exists a map $i:\overline{C(V)} \to X$ as an embedding.
By the construction of $\overline{C(V)}$ in this section, 
if $\overline{C(V)}$ is a complex submanifold in $X$, 
then there exists an embedding $i_{\lambda}: \overline{C_{\lambda}(V)} \to \varphi_{\lambda} (U_{\lambda})$ for each $\lambda \in \Lambda$.

\subsection{Subtorus Actions on Toric Manifolds}
\label{subsec: subtorus actions on toric manifolds}

In this section, we consider a subtorus action on toric manifolds in order to give a Hamiltonian torus action on a complex submanifold given in Section \ref{subsec: construction of toric submanifolds}.

First we define a $k$-dimensional torus action on a complex $n$-dimensional toric manifold $X$.
Given an affine subspace $V= \mathbb{R}p_{1} + \cdots + \mathbb{R}p_{k} + a$ in $\mathbb{R}^{n}$, define a map $i_{V}:T^{k} \to T^{n}$ by 
\begin{equation}
	\label{eq: inclusion of torus}
	i_{V}(t_{1},\ldots, t_{k}) = \left(\prod_{l=1}^{k} t_{l}^{\langle p_{l},e_{1} \rangle}, \ldots, \prod_{l=1}^{k} t_{l}^{\langle p_{l},e_{n} \rangle} \right).
\end{equation}
Recall that the $n$-dimensional torus $T^{n}$ action on $\varphi_{\lambda}(U_{\lambda}) \cong \mathbb{C}^{n}$ is given by 
\begin{center}
$\begin{array}{ccc}
	T^{n} \times \varphi_{\lambda}(U_{\lambda}) &\to& \varphi_{\lambda}(U_{\lambda}) \\
	(t = (t_{1},\ldots, t_{n}), z^{\lambda} = (z^{\lambda}_{1}, \ldots, z^{\lambda}_{n})) &\mapsto& t \cdot z^{\lambda}
\end{array}$
\end{center}
for each $\lambda \in \Lambda$, where $t \cdot z^{\lambda}$ is defined by 
\begin{align*}
	t \cdot z^{\lambda} =& \left( \prod_{j=1}^{n}t_{j}^{Q^{\lambda}_{j1}} z^{\lambda}_{1}, \ldots, \prod_{j=1}^{n}t_{j}^{Q^{\lambda}_{jn}} z^{\lambda}_{n} \right) \\
	=& \left( \prod_{j=1}^{n}t_{j}^{\langle e_{j}, v^{\lambda}_{1} \rangle} z^{\lambda}_{1}, \ldots, \prod_{j=1}^{n}t_{j}^{\langle e_{j}, v^{\lambda}_{n} \rangle} z^{\lambda}_{n} \right).
\end{align*}
This torus action is compatible with the torus action on $(\mathbb{C}^{\ast})^{n}$, which is given by 
\begin{center}
	$\begin{array}{ccc}
		T^{n} \times (\mathbb{C}^{\ast})^{n} &\to& (\mathbb{C}^{\ast})^{n} \\
		((t_{1},\ldots, t_{n}), (z_{1}, \ldots, z_{n})) &\mapsto& (t_{1}z_{1},\ldots, t_{n}z_{n}).
	\end{array}$
	\end{center}
We can describe the $k$-dimensional torus $T^{k}$ action on $X$ by 
\begin{center}
$\begin{array}{ccc}
	T^{k} \times \varphi_{\lambda}(U_{\lambda}) &\to& \varphi_{\lambda}(U_{\lambda}) \\
	(t = (t_{1},\ldots, t_{k}), z^{\lambda}) &\mapsto& i_{V}(t) \cdot z^{\lambda},
\end{array}$
\end{center}
where 
\begin{equation}
	\label{eq: t^k action on C(V)}
	i_{V}(t) \cdot z^{\lambda} 
	= \left( \prod_{l=1}^{k}t_{l}^{\langle p_{l}, v^{\lambda}_{1} \rangle}z^{\lambda}_{1},\ldots, \prod_{l=1}^{k}t_{l}^{\langle p_{l}, v^{\lambda}_{n} \rangle}z^{\lambda}_{n} \right).
\end{equation}

We define the subset $\mathcal{J}_{\lambda} = \{j_{1},\ldots, j_{m}\} \subset \{1,\ldots, n \}$ by
\begin{equation*}
	\mathcal{J}_{\lambda}
	= \{ j \in \{1,\ldots, n \} \mid \langle p_{1},v^{\lambda}_{j} \rangle = \cdots = \langle p_{k},v^{\lambda}_{j} \rangle = 0 \}
\end{equation*}
for $\lambda \in \Lambda$.
If $\mathcal{J}_{\lambda} = \emptyset$, then we may interpret $\mathcal{J}_{\lambda} = \{j_{1},\ldots, j_{m} \}$ as $m=0$.
Note that since the vectors $p_{1},\ldots, p_{k}$ form a basis of the linear part of $V$, we see that 
\begin{equation*}
	\max_{\lambda \in \Lambda} \vert \mathcal{J}_{\lambda} \vert \leq n - k.
\end{equation*}

Since our complex manifold $X$ coincides with the one constructed in \cite{MR4417716} as compact toric manifolds, $X$ is equipped with the torus invariant symplectic form $\omega$ given by Guillemin \cite{MR1293656} (see also \cite[Appendix 2]{MR1301331}).
Moreover, the $T^{n}$-action on $X$ is Hamiltonian with respect to the symplectic form $\omega$. We write the moment map for the Hamiltonian $T^{n}$-action on $X$ is $\mu: X \to (\mathfrak{t}^{n})^{\ast}$.
The fundamental property of moment maps for Hamiltonian torus actions is the convexity theorem \cite{MR642416,MR664117}.

\begin{remark}
	The $T^{k}$-action given in Equation \ref{eq: t^k action on C(V)} is also Hamiltonian with respect to $\omega$ and the moment map for the action is given by $i_{V}^{\ast}\circ \mu: X \to (\mathfrak{t}^{k})^{\ast}$.
\end{remark}
We study the fixed point set of the $T^{k}$-action on $X$.

\begin{lemma}
\label{lemma: fixed points of subtorus action}
	Consider the $T^{k}$-action on $X$ defined above.
	The fixed point set of the $T^{k}$-action on $\varphi_{\lambda}(U_{\lambda})$ is $\{z^{\lambda} \in \varphi_{\lambda}(U_{\lambda}) \mid z^{\lambda}_{i} = 0, i \notin \mathcal{J}_{\lambda} \}$.
\begin{proof}
For simplicity, we give the proof for the case when $\mathcal{J}_{\lambda} = \{i \}$.

If $\langle p_{l},v^{\lambda}_{i} \rangle = 0$ for $l = 1,\ldots, k$, then we have 
\begin{align*}
	i_{V}(t) \cdot (0, \ldots, 0, z^{\lambda}_{i}, 0, 0\ldots,0)
	&= 
	\left( 0, \ldots, 0, \prod_{l=1}^{k}t_{l}^{\langle p_{l},v^{\lambda}_{i} \rangle} z^{\lambda}_{i}, 0,\ldots, 0 \right) \\
	&= 
	(0, \ldots, 0, z^{\lambda}_{i}, 0,\ldots, 0).
\end{align*}
Thus, the point $(0, \ldots, 0, z^{\lambda}_{i}, 0, \ldots,0)$ is a fixed point of the $T^{k}$-action on $\varphi_{\lambda}(U_{\lambda})$.

Conversely, if $(0, \ldots, 0, z^{\lambda}_{i}, 0, \ldots,0)$ is a fixed point of the $T^{k}$-action on $\varphi_{\lambda}(U_{\lambda})$, then we have 
\begin{equation*}
	\left( 0, \ldots, 0, \prod_{l=1}^{k}t_{l}^{\langle p_{l},v^{\lambda}_{i} \rangle} z^{\lambda}_{i}, 0,\ldots, 0 \right)
	= 
	(0, \ldots, 0, z^{\lambda}_{i}, 0,\ldots, 0).
\end{equation*}
Since $(t_{1},\ldots, t_{k}) \in T^{k}$, we see that $\langle p_{l},v^{\lambda}_{i} \rangle = 0$ for $l = 1,\ldots, k$.

Note that if $\mathcal{J}_{\lambda} = \emptyset$, we see that the fixed point of the $T^{k}$-action on $\varphi_{\lambda}(U_{\lambda})$ is $(0,\ldots,0) \in \varphi_{\lambda}(U_{\lambda})$.
\end{proof}
\end{lemma}

Note that 
the set $\{z^{\lambda} \in \varphi_{\lambda}(U_{\lambda}) \mid z^{\lambda}_{i} = 0, i \notin \mathcal{J}_{\lambda} \}$ corresponds to the $m$-face defined by the direction vectors $v^{\lambda}_{j_{1}},\ldots, v^{\lambda}_{j_{m}}$ in $\Delta$ for $j_{1}, \ldots, j_{m} \in \mathcal{J}_{\lambda}$.
If $\mathcal{J}_{\lambda} = \emptyset$, then the set $\{z^{\lambda} \in \varphi_{\lambda}(U_{\lambda}) \mid z^{\lambda}_{i} = 0, i \notin \mathcal{J}_{\lambda} \} = \{(0,\ldots,0) \in \varphi_{\lambda
}(U_{\lambda}) \}$ corresponds to the vertex $\lambda$, which is a $0$-face in $\Delta$.
In particular, $\{z^{\lambda} \in \varphi_{\lambda}(U_{\lambda}) \mid z^{\lambda}_{i} = 0, i \notin \mathcal{J}_{\lambda} \} \neq \emptyset$ for any $\mathcal{J}_{\lambda}$.

\subsection{Torus Actions on Torus-equivariantly Embedded Toric Manifolds}
\label{subsec: torus actions on toric submanifolds}

After giving torus actions on $\overline{C(V)}$, we consider the image of the moment map for the torus action.

Under the following diagram;
\begin{equation*}
	\begin{tikzcd}
		T^{n} & \times & \varphi_{\lambda}(U_{\lambda}) \arrow[r] & \varphi_{\lambda}(U_{\lambda}) \\
		T^{k} \arrow[u, "i_{V}"'] & \times & \overline{C_{\lambda}(V)} \arrow[r] \arrow[u, "i_{\lambda}"'] & \overline{C_{\lambda}(V)} \arrow[u, "i_{\lambda}"'],
	\end{tikzcd}
\end{equation*}
a $T^{k}$-action on $\overline{C(V)}$ is defined by 
\begin{center}
	$\begin{array}{ccc}
		T^{k} \times \overline{C_{\lambda}(V)} &\to& \overline{C_{\lambda}(V)} \\
		(t = (t_{1},\ldots, t_{k}), z^{\lambda}) &\mapsto& i_{V}(t)\cdot z^{\lambda},
	\end{array}$
\end{center}
which makes the above diagram commutative.

From Lemma \ref{lemma: fixed points of subtorus action} and the definition of $\overline{C_{\lambda}(V)}$, the following lemma is obvious.

\begin{lemma}
	\label{lemma: T^k action on toric submanifolds}
	The fixed point set of the $T^{k}$-action on $\overline{C_{\lambda}(V)}$ is $\{z^{\lambda} \in \varphi_{\lambda}(U_{\lambda}) \mid z^{\lambda}_{i} = 0, i \notin \mathcal{J}_{\lambda} \} \cap \overline{C_{\lambda}(V)}$.
\end{lemma}
It is clear that if $\{z^{\lambda} \in \varphi_{\lambda}(U_{\lambda}) \mid z^{\lambda}_{i} = 0, i \notin \mathcal{J}_{\lambda} \} \cap \overline{C_{\lambda}(V)} \neq \emptyset$, then there exists a fixed point of the $T^{k}$-action on $\overline{C_{\lambda}(V)}$.

\begin{lemma}
	Assume that $k \geq 1$.
	If there exists $j = k+1,\ldots, n$ such that $\langle u^{\lambda}_{i},q_{j} \rangle > 0$ (or $\langle u^{\lambda}_{i},q_{j} \rangle < 0$) for all $i=1,\ldots,n$, then $\{z^{\lambda} \in \varphi_{\lambda}(U_{\lambda}) \mid z^{\lambda}_{i} = 0, i \notin \mathcal{J}_{\lambda} \} \cap \overline{C_{\lambda}(V)} = \emptyset$, i.e., there is no fixed point of the $T^{k}$-action on $\overline{C_{\lambda}(V)}$.
\begin{proof}
	For simplicity, we assume that there exists $j=k+1,\ldots,n$ such that $\langle u^{\lambda}_{i},q_{j} \rangle > 0$ for all $i=1,\ldots,n$.
	This implies that $\mathcal{I}^{+}_{\lambda,j} = \{1,\ldots, n\}$ and $\mathcal{I}^{-}_{\lambda,j} = \emptyset$.
	As we noted in the definition of $f^{\lambda}_{j}(z^{\lambda})$, we obtain 
	\begin{align*}
		f^{\lambda}_{j}(z^{\lambda})
		&=
		\prod_{i \in \mathcal{I}^{+}_{\lambda,j}}(z^{\lambda}_{i})^{\langle u^{\lambda}_{i},q_{j} \rangle} 
		- e^{\langle a,q_{j} \rangle} \prod_{i \in \mathcal{I}^{-}_{\lambda,j}}(z^{\lambda}_{i})^{-\langle u^{\lambda}_{i},q_{j} \rangle} \\
		&=
		\prod_{i=1}^{n}(z^{\lambda}_{i})^{\langle u^{\lambda}_{i},q_{j} \rangle} - e^{\langle a,q_{j} \rangle}.
	\end{align*}
	Since $e^{\langle a,q_{j} \rangle} \neq 0$, $f^{\lambda}_{j}(z^{\lambda}) = 0$ implies that $z^{\lambda}_{1}z^{\lambda}_{2} \cdots z^{\lambda}_{n} \neq 0$, i.e.,
	\begin{equation*}
		\{f^{\lambda}_{j}(z^{\lambda}) = 0 \} \subset \{z^{\lambda}_{1}z^{\lambda}_{2}\cdots z^{\lambda}_{n} \neq 0 \}.
	\end{equation*} 

	Recall that $\vert \mathcal{J}_{\lambda} \vert \leq n-k$ for any $\lambda \in \Lambda$. 
	If $k \geq 1$, then there exists $i_{0} \not\in \mathcal{J}_{\lambda}$ such that 
	\begin{equation*}
		\{z^{\lambda} \in \varphi_{\lambda}(U_{\lambda}) \mid z^{\lambda}_{i} = 0, i\not\in \mathcal{J}_{\lambda} \} \subset \{z^{\lambda}_{i_{0}} = 0 \}.
	\end{equation*}

	It is clear that $\{z^{\lambda}_{1}z^{\lambda}_{2}\cdots z^{\lambda}_{n} \neq 0 \} \cap \{z^{\lambda}_{i_{0}} = 0 \} = \emptyset$, which implies that $\{z^{\lambda} \in \varphi_{\lambda}(U_{\lambda}) \mid z^{\lambda}_{i} = 0, i \notin \mathcal{J}_{\lambda} \} \cap \overline{C_{\lambda}(V)} = \emptyset$.
	Since $\overline{C_{\lambda}(V)} = \bigcap_{j=k+1}^{n}\{f^{\lambda}_{j}(z^{\lambda}) = 0 \}$, we obtain the desired result.
\end{proof}
\end{lemma}

Note that if $\langle u^{\lambda}_{i},q_{j} \rangle > 0$ (or $\langle u^{\lambda}_{i},q_{j} \rangle < 0$) for all $i=1,\ldots,n$ and some $j=k+1,\ldots,n$, then $\overline{C_{\lambda}(V)} \subset \varphi_{\lambda}(\check{U}_{\lambda})$.

\begin{remark}
	If $\overline{C(V)}$ is a complex submanifold in $X$, the $T^{k}$-action on $\overline{C(V)}$ is actually Hamiltonian with respect to the symplectic form $i^{\ast}\omega$ on $\overline{C(V)}$ for the inclusion $i: \overline{C(V)} \to X$.
If $\mu:X \to (\mathfrak{t}^{n})^{\ast}$ is the moment map for the $T^{n}$-action on $X$, then we can obtain the moment map for the $T^{k}$-action on $\overline{C(V)}$ by $\overline{\mu} = i_{V}^{\ast} \circ \mu \circ i: \overline{C(V)} \to (\mathfrak{t}^{k})^{\ast}$.
\end{remark}

We further examine the fixed points of the $T^{k}$-action on $\overline{C(V)}$.

Since the map $i_{V}:T^{k} \to T^{n}$ is defined by (\ref{eq: inclusion of torus}), we can write the pull back $i_{V}^{\ast}: (\mathfrak{t}^{n})^{\ast} \to (\mathfrak{t}^{k})^{\ast}$ as 
	\begin{equation}
		\label{eq: i^*_{V}}
		i_{V}^{\ast} (\xi)= \left(\langle p_{1},\xi \rangle, \ldots, \langle p_{k},\xi \rangle \right).
	\end{equation}
Note that this map $i_{V}^{\ast}$ is a surjective linear map.
Hence, we have 
\begin{equation*}
	i_{V}^{\ast}(\mu(X))= \{(\langle p_{1},\xi \rangle, \ldots, \langle p_{k},\xi\rangle) \mid \xi \in \mu(X) \} \subset (\mathfrak{t}^{k})^{\ast}.
\end{equation*}
Since $X$ is a toric manifold, $\mu (X) = \Delta$ is a Delzant polytope. In this situation, we obtain the followings:
\begin{cor}
	\label{cor: delzant to convex}
	Let $\Delta$ be a Delzant polytope and $i_{V}^{\ast}: (\mathfrak{t}^{n})^{\ast} \to (\mathfrak{t}^{k})^{\ast}$ be the map defined in Equation \ref{eq: i^*_{V}}.
	Then, $i_{V}^{\ast}(\Delta)$ is a convex polytope.
\begin{proof}
	As we noted, the map $i_{V}^{\ast}$ is a linear map.
	By the definition of Delzant polytopes, $\Delta$ is a convex polytope. 
	By Lemma \ref{lemma: convex polytope linear map}, $i_{V}^{\ast}(\Delta)$ is a convex polytope.
\end{proof}
\end{cor}

\begin{lemma}
	\label{lemma: fixed point in face}
	Suppose that $\mathcal{J}_{\lambda} = \{ j_{1},\ldots, j_{m}\}$.
	Let $\mathcal{F}_{\lambda}$ be an $m$-face of $\Delta$ defined by the direction vectors $v^{\lambda}_{j_{1}}, \ldots, v^{\lambda}_{j_{m}}$.
	Then, $i_{V}^{\ast}(\xi) = i_{V}^{\ast}(\lambda)$ holds for any $\xi \in \mathcal{F}_{\lambda}$.

\begin{proof}
	Since $\xi, \lambda \in \mathcal{F}_{\lambda}$, we have 
	\begin{equation*}
		\xi - \lambda = \sum_{l=1}^{m}\alpha_{l}v^{\lambda}_{j_{l}}
	\end{equation*}
	for some $\alpha_{1}, \ldots, \alpha_{m} \in \mathbb{R}$. We calculate 
	\begin{align*}
		i_{V}^{\ast}(\xi) - i_{V}^{\ast}(\lambda) 
		&= 
		i_{V}^{\ast}(\xi - \lambda) \\
		&= 
		\left(\langle p_{1}, \sum_{l=1}^{m}\alpha_{l}v^{\lambda}_{j_{l}} \rangle, \ldots, \langle p_{k}, \sum_{l=1}^{m}\alpha_{l}v^{\lambda}_{j_{l}} \rangle \right) \\
		&=
		\left( \sum_{l=1}^{m}\alpha_{l}\langle p_{1},v^{\lambda}_{j_{l}} \rangle, \ldots, \sum_{l=1}^{m}\alpha_{l}\langle p_{k},v^{\lambda}_{j_{l}} \rangle \right) \\
		&= (0,\ldots, 0).
	\end{align*}
	Thus we obtain $i_{V}^{\ast}(\xi) = i_{V}^{\ast}(\lambda)$ for any $\xi \in \mathcal{F}_{\lambda}$.
\end{proof}
\end{lemma}

By the definition of Delzant polytopes,
we can take the direction vectors $v^{\lambda}_{1},\ldots, v^{\lambda}_{n} \in \mathbb{Z}^{n}$ from the vertex $\lambda$ of a Delzant polytope $\Delta$ and the vectors $v^{\lambda}_{1},\ldots, v^{\lambda}_{n} \in \mathbb{Z}^{n}$ can be chosen as a basis of $\mathbb{Z}^{n}$.
We define the cone $\mathcal{C}_{\lambda}$ by 
\begin{equation*}
	\mathcal{C}_{\lambda}
	=
	\mathbb{R}_{\geq 0}v^{\lambda}_{1}
	+ \cdots +
	\mathbb{R}_{\geq 0}v^{\lambda}_{n}
	\subset (\mathfrak{t}^{n})^{\ast} \cong \mathbb{R}^{n}.
\end{equation*}
In other words, $\mathcal{C}_{\lambda}$ is generated by $\{v^{\lambda}_{1},\ldots,v^{\lambda}_{n} \}$.
Since the map $i_{V}^{\ast}$ is linear, by Lemma \ref{lemma: cone linear map}, $i_{V}^{\ast}(\mathcal{C}_{\lambda})$ is the cone generated by $\{i_{V}^{\ast}(v^{\lambda}_{1}),\ldots, i_{V}^{\ast}(v^{\lambda}_{n}) \}$. The cone $i_{V}^{\ast}(\mathcal{C}_{\lambda})$ can be written concretely by 
\begin{equation*}
	i_{V}^{\ast}(\mathcal{C}_{\lambda})
	=
	\mathbb{R}_{\geq 0}i_{V}^{\ast}(v^{\lambda}_{1}) + \cdots + \mathbb{R}_{\geq 0}i_{V}^{\ast}(v^{\lambda}_{n}) \subset (\mathfrak{t}^{k})^{\ast} \cong \mathbb{R}^{k}.
\end{equation*}
\begin{defin}
	Let $\Delta$ be a Delzant polytope and $\lambda$ be a vertex in $\Delta$.
	The point $i_{V}^{\ast}(\lambda)$ is a vertex in the convex polytope $i_{V}^{\ast}(\Delta)$ if $0 = (0,\ldots,0) \in i_{V}^{\ast}(\mathcal{C}_{\lambda}) \subset (\mathfrak{t}^{k})^{\ast}$ is a vertex in the sense of Definition \ref{def: vertex of cone}.
\end{defin}

We have the relation between the vectors $p_{1},\ldots,p_{k}$, $q_{k+1}, \ldots, q_{n}$, $v^{\lambda}_{1},\ldots,v^{\lambda}_{n}$, and $u^{\lambda}_{1},\ldots, u^{\lambda}_{n}$.
\begin{lemma}
	\label{lemma: p,v,u,q}
	For any $l=1,\ldots,k$ and any $j=k+1,\ldots, n$, 
	\begin{equation}
		\label{eq: p,v,u,q}
		\sum_{i=1}^{n} \langle p_{l},v^{\lambda}_{i} \rangle \langle u^{\lambda}_{i}, q_{j} \rangle = 0
	\end{equation}
	holds.
\begin{proof}
	Instead of Equation \ref{eq: p,v,u,q}, we show the matrix equation 
	\begin{equation}
		\label{eq: p,v,u,q matrix}
		\left[
			\begin{matrix}
				{^t p_{1}} \\
				\vdots \\
				{^t p_{k}}
			\end{matrix}
		\right]
		\left[
			\begin{matrix}
				v^{\lambda}_{1} & \cdots & v^{\lambda}_{n}
			\end{matrix}
		\right]
		\left[
			\begin{matrix}
				{^t u^{\lambda}_{1}} \\
				\vdots \\
				{^t u^{\lambda}_{n}}
			\end{matrix}
		\right]
		\left[
			\begin{matrix}
				q_{k+1} & \cdots & q_{n}
			\end{matrix}
		\right]
		= 0.
	\end{equation}
	
	From Lemma \ref{lemma: direction vectors vs normal vectors}, we have obtained 
	\begin{equation}
		\label{eq: v,u}
		\left[
			\begin{matrix}
				v^{\lambda}_{1} & \cdots & v^{\lambda}_{n}
			\end{matrix}
		\right]
		\left[
			\begin{matrix}
				{^t u^{\lambda}_{1}} \\
				\vdots \\
				{^t u^{\lambda}_{n}}
			\end{matrix}
		\right]
		= E_{n}.
	\end{equation}
	Since the vectors $q_{k+1},\ldots, q_{n}$ are taken to be an orthogonal basis of the orthogonal subspace to $V$, we have 
	\begin{equation}
		\label{eq: p,q}
		\left[
			\begin{matrix}
				{^t p_{1}} \\
				\vdots \\
				{^t p_{k}}
			\end{matrix}
		\right]
		\left[
			\begin{matrix}
				q_{k+1} & \cdots & q_{n}
			\end{matrix}
		\right] = 0.
	\end{equation}
	From Equation \ref{eq: v,u} and Equation \ref{eq: p,q}, we calculate 
	\begin{align*}
		\left[
			\begin{matrix}
				{^t p_{1}} \\
				\vdots \\
				{^t p_{k}}
			\end{matrix}
		\right]
		\left[
			\begin{matrix}
				v^{\lambda}_{1} & \cdots & v^{\lambda}_{n}
			\end{matrix}
		\right]
		\left[
			\begin{matrix}
				{^t u^{\lambda}_{1}} \\
				\vdots \\
				{^t u^{\lambda}_{n}}
			\end{matrix}
		\right]
		\left[
			\begin{matrix}
				q_{k+1} & \cdots & q_{n}
			\end{matrix}
		\right]
		&=
		\left[
			\begin{matrix}
				{^t p_{1}} \\
				\vdots \\
				{^t p_{k}}
			\end{matrix}
		\right]
		\left[
			\begin{matrix}
				q_{k+1} & \cdots & q_{n}
			\end{matrix}
		\right] \\
		&= 0.
	\end{align*}
	Hence, we obtain Equation \ref{eq: p,v,u,q matrix}.
\end{proof}
\end{lemma}

From Lemma \ref{lemma: p,v,u,q}, we see that 
\begin{equation*}
	\left[
		\begin{matrix}
			\displaystyle
			\sum_{i=1}^{n}\langle p_{1},v^{\lambda}_{i} \rangle \langle u^{\lambda}_{i}, q_{j} \rangle \\
			\vdots \\
			\displaystyle
			\sum_{i=1}^{n}\langle p_{k},v^{\lambda}_{i} \rangle \langle u^{\lambda}_{i}, q_{j} \rangle
		\end{matrix}
	\right]
	=0
\end{equation*}
holds for any $j=k+1,\ldots,n$.
This is equivalent to the equation:
\begin{equation*}
	\sum_{i=1}^{n}
	\langle u^{\lambda}_{i}, q_{j} \rangle
	\left[
		\begin{matrix}
			\displaystyle
			\langle p_{1},v^{\lambda}_{i} \rangle \\
			\vdots \\
			\displaystyle
			\langle p_{k},v^{\lambda}_{i} \rangle 
		\end{matrix}
	\right]
	=0.
\end{equation*}
From (\ref{eq: i^*_{V}}), the following equation 
\begin{equation}
	\label{eq: i,v,u,q}
	\sum_{i=1}^{n}\langle u^{\lambda}_{i}, q_{j} \rangle i_{V}^{\ast}(v^{\lambda}_{i}) = 0
\end{equation}
holds for any $j=k+1,\ldots,n$.
Since the set $\mathcal{J}_{\lambda} \subset \{1,\ldots, n \}$ was defined by 
\begin{equation*}
	\mathcal{J}_{\lambda}
	=
	\{
		i \in \{1,\ldots, n \} 
		\mid 
		\langle p_{1},v^{\lambda}_{i} \rangle = \cdots = \langle p_{k},v^{\lambda}_{i} \rangle = 0
	\},
\end{equation*}
Equation \ref{eq: i,v,u,q} can be written as 
\begin{equation}
	\label{eq: i,v,u,q refinement}
	\sum_{i \not\in \mathcal{J}_{\lambda}} \langle u^{\lambda}_{i}, q_{j} \rangle i_{V}^{\ast}(v^{\lambda}_{i}) = 0.
\end{equation}
Note that since $\vert \mathcal{J}_{\lambda} \vert \leq n-k$, the number of the terms in the left hand side of Equation \ref{eq: i,v,u,q refinement} should be greater than or equal to $k$.

\begin{lemma}
	\label{lemma: star conditions}
	Fix $j = k+1,\ldots,n$.
	If the point $i_{V}^{\ast}(\lambda)$ is a vertex in the convex polytope $i_{V}^{\ast}(\Delta)$, then $\{\langle u^{\lambda}_{i},q_{j} \rangle\}_{i \not\in \mathcal{J}_{\lambda}}$ satisfies either of the following conditions:
	\begin{enumerate} 
		\item $\langle u^{\lambda}_{i},q_{j} \rangle = 0$ holds for any $i \not\in \mathcal{J}_{\lambda}$,
		\item there exist $i_{j} \neq i^{\prime}_{j} \not\in \mathcal{J}_{\lambda}$ such that $\langle u^{\lambda}_{i_{j}},q_{j} \rangle \langle u^{\lambda}_{i^{\prime}_{j}},q_{j} \rangle < 0$.
	\end{enumerate}
\begin{proof}
	Regarding Equation \ref{eq: i,v,u,q refinement} as a linear combination of the vectors $i_{V}^{\ast}(v^{\lambda}_{i}) \; (i \not\in \mathcal{J}_{\lambda})$, the coefficients $\langle u^{\lambda}_{i},q_{j} \rangle\; (i \not\in \mathcal{J}_{\lambda})$ satisfy at least one of the following cases:
	\begin{itemize}
		\item $\langle u^{\lambda}_{i},q_{j} \rangle \geq 0$ holds for any $i \not\in \mathcal{J}_{\lambda}$,
		\item $\langle u^{\lambda}_{i},q_{j} \rangle \leq 0$ holds for any $i \not\in \mathcal{J}_{\lambda}$,
		\item there exist $i_{j} \neq i^{\prime}_{j} \not\in \mathcal{J}_{\lambda}$ such that $\langle u^{\lambda}_{i_{j}},q_{j} \rangle \langle u^{\lambda}_{i^{\prime}_{j}},q_{j} \rangle < 0$.
	\end{itemize}
	If the point $i_{V}^{\ast}(\lambda)$ is a vertex in the convex polytope $i_{V}^{\ast}(\Delta)$, then by Definition \ref{def: vertex of cone} and Equation \ref{eq: i,v,u,q refinement} the first and second cases can be written as $\langle u^{\lambda}_{i},q_{j} \rangle = 0$ holds for any $i \not\in \mathcal{J}_{\lambda}$.
\end{proof}
\end{lemma}

To check whether the fixed point set of the $T^{k}$-action on $\overline{C_{\lambda}(V)}$ is empty or not, we can use Lemma \ref{lemma: star conditions}.

\begin{lemma}
	\label{lemma: clubsuit 2}
	Fix $j=k+1,\ldots,n$. 
	If the point $i_{V}^{\ast}(\lambda)$ is a vertex in the convex polytope $i_{V}^{\ast}(\Delta)$ 
	and if there exist $i_{j} \neq i^{\prime}_{j}$ such that $\langle u^{\lambda}_{i_{j}},q_{j} \rangle \langle u^{\lambda}_{i^{\prime}_{j}},q_{j} \rangle < 0$,
	then  
	\begin{equation}
		\label{eq: clubsuit 2}
		\{z^{\lambda} \in \varphi_{\lambda}(U_{\lambda})\mid z^{\lambda}_{i_{j}} = z^{\lambda}_{i^{\prime}_{j}} = 0 \} 
		\subset 
		\{z^{\lambda} \in \varphi_{\lambda}(U_{\lambda})\mid f^{\lambda}_{j}(z^{\lambda}) = 0 \} 
	\end{equation}
	holds.
\begin{proof}
	For simplicity, we assume that $\langle u^{\lambda}_{i_{j}},q_{j} \rangle >0$, $\langle u^{\lambda}_{i^{\prime}_{j}},q_{j} \rangle < 0$.
	This implies that $i_{j} \in \mathcal{I}^{+}_{\lambda,j} \setminus \mathcal{I}^{0}_{\lambda,j}$ and $i^{\prime}_{j} \in \mathcal{I}^{-}_{\lambda,j} \setminus \mathcal{I}^{0}_{\lambda,j}$.
	From Equation \ref{eq: defining eqautions for toric submanifolds}, we obtain Equation \ref{eq: clubsuit 2}. 
\end{proof}
\end{lemma}

As a corollary to Lemma \ref{lemma: clubsuit 2}, if the point $i_{V}^{\ast}(\lambda)$ is a vertex in the convex polytope $i_{V}^{\ast}(\Delta)$ and if $\{\langle u^{\lambda}_{i},q_{j} \rangle\}_{i \not\in \mathcal{J}_{\lambda}}$ satisfies the condition $(2)$ in Lemma \ref{lemma: star conditions} for a fixed $j$, then we obtain 
\begin{equation*}
	\{z^{\lambda} \in \varphi_{\lambda}(U_{\lambda})\mid z^{\lambda}_{i_{j}} = z^{\lambda}_{i^{\prime}_{j}} = 0, i_{j} \neq i^{\prime}_{j} \not\in \mathcal{J}_{\lambda} \} 
	\subset 
	\{z^{\lambda} \in \varphi_{\lambda}(U_{\lambda})\mid f^{\lambda}_{j}(z^{\lambda}) = 0 \} 
\end{equation*}
for $i_{j} \neq i^{\prime}_{j} \not\in \mathcal{J}_{\lambda}$ appearing in the statement of the condition $(2)$ in Lemma \ref{lemma: star conditions}.

When $\{\langle u^{\lambda}_{i},q_{j} \rangle\}_{i \not\in \mathcal{J}_{\lambda}}$ satisfies the condition $(1)$ in Lemma \ref{lemma: star conditions}, we obtain the following result:

\begin{prop}
	\label{prop: clubsuit 1}
	Assume that the point $i_{V}^{\ast}(\lambda)$ is a vertex in the convex polytope $i_{V}^{\ast}(\Delta)$.
	We define a point $\tilde{z}^{\lambda}=(\tilde{z}^{\lambda}_{1}, \ldots, \tilde{z}^{\lambda}_{n}) \in \mathbb{C}^{n}$ by setting 
	\begin{equation*}
		\tilde{z}^{\lambda}_{i}
		=
		\begin{cases}
			e^{\langle a,v^{\lambda}_{i} \rangle}, & i \in \mathcal{J}_{\lambda}, \\
			0, & i \not\in \mathcal{J}_{\lambda}.
		\end{cases}
	\end{equation*}

	If there exists $j_{0} = k+1,\ldots,n$ such that $\{\langle u^{\lambda}_{i},q_{j_{0}} \rangle\}_{i \not\in \mathcal{J}_{\lambda}}$ satisfies the condition $(1)$ in Lemma \ref{lemma: star conditions}, 
	then we obtain 
	\begin{equation*}
		\tilde{z}^{\lambda} \in \{z^{\lambda} \in \varphi_{\lambda}(U_{\lambda}) \mid z^{\lambda}_{i} = 0, i\not\in \mathcal{J}_{\lambda} \} \cap \overline{C_{\lambda}(V)},
	\end{equation*}
	i.e., the set $\{z^{\lambda} \in \varphi_{\lambda}(U_{\lambda}) \mid z^{\lambda}_{i} = 0, i\not\in \mathcal{J}_{\lambda} \} \cap \overline{C_{\lambda}(V)}$ is not empty. 
	In particular, there exists a fixed point of the $T^{k}$-action on $\overline{C_{\lambda}(V)}$.
\begin{proof}
	It is clear that $\tilde{z}^{\lambda} \in \{z^{\lambda} \in \varphi_{\lambda}(U_{\lambda}) \mid z^{\lambda}_{i} = 0, i\not\in \mathcal{J}_{\lambda} \}$.
	We show that $\tilde{z}^{\lambda} \in \overline{C_{\lambda}(V)} = \bigcap_{j=k+1}^{n}\{z^{\lambda}\in \varphi_{\lambda}(U_{\lambda}) \mid f^{\lambda}_{j}(z^{\lambda})=0 \}$.

	Since we assume that the point $i_{V}^{\ast}(\lambda)$ is a vertex in the convex polytope $i_{V}^{\ast}(\Delta)$, $\{ \langle u^{\lambda}_{i},q_{j} \rangle \}_{i \not\in \mathcal{J}_{\lambda}}$ satisfies either of the condition $(1)$ or the condition $(2)$ in Lemma \ref{lemma: star conditions} for each fixed $j$.
	We use the result to check that $\tilde{z}^{\lambda} \in \{z^{\lambda}\in \varphi_{\lambda}(U_{\lambda}) \mid f^{\lambda}_{j}(z^{\lambda})=0 \}$ for any $j$.

	If $\{ \langle u^{\lambda}_{i},q_{j} \rangle \}_{i \not\in \mathcal{J}_{\lambda}}$ satisfies the condition $(1)$ in Lemma \ref{lemma: star conditions} for some $j$, then 
	$i \not\in \mathcal{J}_{\lambda}$ implies $\langle u^{\lambda}_{i},q_{j} \rangle = 0$, i.e., $i \in \mathcal{I}^{0}_{\lambda,j}$ for such $j$.
	By considering the contraposition, $\mathcal{I}^{+}_{\lambda,j} \cup \mathcal{I}^{-}_{\lambda,j} \setminus \mathcal{I}^{0}_{\lambda,j} \subset \mathcal{J}_{\lambda}$.
	We calculate $\prod_{i \in \mathcal{I}^{+}_{\lambda,j}} (\tilde{z}^{\lambda}_{i})^{\langle u^{\lambda}_{i},q_{j} \rangle}$ and $\prod_{i \in \mathcal{I}^{-}_{\lambda,j}} (\tilde{z}^{\lambda}_{i})^{-\langle u^{\lambda}_{i},q_{j} \rangle}$ as 
	\begin{align*}
		&\prod_{i \in \mathcal{I}^{+}_{\lambda,j}} (\tilde{z}^{\lambda}_{i})^{\langle u^{\lambda}_{i},q_{j} \rangle} 
		=
		\prod_{i \in \mathcal{I}^{+}_{\lambda,j} \setminus \mathcal{I}^{0}_{\lambda,j}} (\tilde{z}^{\lambda}_{i})^{\langle u^{\lambda}_{i},q_{j} \rangle} 
		=
		\prod_{i \in \mathcal{I}^{+}_{\lambda,j} \setminus \mathcal{I}^{0}_{\lambda,j}} (e^{\langle a,v^{\lambda}_{i}\rangle})^{\langle u^{\lambda}_{i},q_{j} \rangle}, \\
		&\prod_{i \in \mathcal{I}^{-}_{\lambda,j}} (\tilde{z}^{\lambda}_{i})^{-\langle u^{\lambda}_{i},q_{j} \rangle}  
		=
		\prod_{i \in \mathcal{I}^{-}_{\lambda,j} \setminus \mathcal{I}^{0}_{\lambda,j}} (\tilde{z}^{\lambda}_{i})^{-\langle u^{\lambda}_{i},q_{j} \rangle}  
		= 
		\prod_{i \in \mathcal{I}^{-}_{\lambda,j} \setminus \mathcal{I}^{0}_{\lambda,j}} (e^{\langle a,v^{\lambda}_{i} \rangle})^{-\langle u^{\lambda}_{i},q_{j} \rangle}
		\neq 0.
	\end{align*}
	Moreover, since $i \in \mathcal{I}^{0}_{\lambda,j}$ means that $\langle u^{\lambda}_{i},q_{j} \rangle = 0$, we obtain 
	\begin{align*}
		&\prod_{i \in \mathcal{I}^{-}_{\lambda,j} \setminus \mathcal{I}^{0}_{\lambda,j}} (e^{\langle a,v^{\lambda}_{i} \rangle})^{-\langle u^{\lambda}_{i},q_{j} \rangle}
		=
		\prod_{i \in \mathcal{I}^{-}_{\lambda,j}} (e^{\langle a,v^{\lambda}_{i} \rangle})^{-\langle u^{\lambda}_{i},q_{j} \rangle}, \\
		&\prod_{i \in \mathcal{I}^{+}_{\lambda,j}\setminus \mathcal{I}^{0}_{\lambda,j}} (e^{\langle a,v^{\lambda}_{i} \rangle})^{\langle u^{\lambda}_{i},q_{j} \rangle}
		= 
		\prod_{i \in \mathcal{I}^{+}_{\lambda,j}} (e^{\langle a,v^{\lambda}_{i} \rangle})^{\langle u^{\lambda}_{i},q_{j} \rangle}.
	\end{align*} 
	As we noted that $\mathcal{I}^{+}_{\lambda,j} \cup \mathcal{I}^{-}_{\lambda,j} = \{1,\ldots, n \}$, we can calculate 
	\begin{align*}
		\frac{\prod_{i \in \mathcal{I}^{+}_{\lambda,j}} (\tilde{z}^{\lambda}_{i})^{\langle u^{\lambda}_{i},q_{j} \rangle}}{\prod_{i \in \mathcal{I}^{-}_{\lambda,j}} (\tilde{z}^{\lambda}_{i})^{-\langle u^{\lambda}_{i},q_{j} \rangle}}
		&=
		\frac{\prod_{i \in \mathcal{I}^{+}_{\lambda,j} \setminus \mathcal{I}^{0}_{\lambda,j}} (e^{\langle a, v^{\lambda}_{i} \rangle})^{\langle u^{\lambda}_{i},q_{j} \rangle}}{\prod_{i \in \mathcal{I}^{-}_{\lambda,j} \setminus \mathcal{I}^{0}_{\lambda,j}} (e^{\langle a, v^{\lambda}_{i} \rangle})^{-\langle u^{\lambda}_{i},q_{j} \rangle}} \\
		&=
		\frac{\prod_{i \in \mathcal{I}^{+}_{\lambda,j}} (e^{\langle a, v^{\lambda}_{i} \rangle})^{\langle u^{\lambda}_{i},q_{j} \rangle}}{\prod_{i \in \mathcal{I}^{-}_{\lambda,j}} (e^{\langle a, v^{\lambda}_{i} \rangle})^{-\langle u^{\lambda}_{i},q_{j} \rangle}} \\
		&=
		\prod_{i \in \mathcal{I}^{+}_{\lambda,j} \cup \mathcal{I}^{-}_{\lambda,j}} (e^{\langle a, v^{\lambda}_{i} \rangle})^{\langle u^{\lambda}_{i},q_{j} \rangle} \\ 
		&=
		\prod_{i =1}^{n} (e^{\langle a,v^{\lambda}_{i} \rangle})^{\langle u^{\lambda}_{i},q_{j} \rangle} \\
		&=
		e^{\langle a,q_{j} \rangle}.
	\end{align*} 
	This calculation implies that $f^{\lambda}_{j}(\tilde{z}^{\lambda}) = 0$ for $j$ such that $\{ \langle u^{\lambda}_{i},q_{j} \rangle \}_{i \not\in \mathcal{J}_{\lambda}}$ satisfies the condition $(1)$ in Lemma \ref{lemma: star conditions}.

	If $\{ \langle u^{\lambda}_{i},q_{j} \rangle \}_{i \not\in \mathcal{J}_{\lambda}}$ does not satisfy the condition $(1)$ in Lemma \ref{lemma: star conditions} for some $j$, i.e., if $\{ \langle u^{\lambda}_{i},q_{j} \rangle \}_{i \not\in \mathcal{J}_{\lambda}}$ satisfies the condition $(2)$ in Lemma \ref{lemma: star conditions} for some $j$, 
	then by Lemma \ref{lemma: clubsuit 2}, there exist $i_{j} \neq i^{\prime}_{j} \not\in \mathcal{J}_{\lambda}$ such that 
	Equation \ref{eq: clubsuit 2} holds for such $j$.
	By the definition of $\tilde{z}^{\lambda}$, $\tilde{z}^{\lambda} \in \{z^{\lambda} \in \varphi_{\lambda}(U_{\lambda}) \mid z^{\lambda}_{i_{j}} = z^{\lambda}_{i^{\prime}_{i}} = 0 \}$ holds for such $j$, which implies that $\tilde{z}^{\lambda} \in \{z^{\lambda}\in \varphi_{\lambda}(U_{\lambda}) \mid f^{\lambda}_{j}(z^{\lambda})=0 \}$.
	
	From the above discussion, $\tilde{z}^{\lambda} \in \{z^{\lambda}\in \varphi_{\lambda}(U_{\lambda}) \mid f^{\lambda}_{j}(z^{\lambda})=0 \}$ holds for any $j = k+1,\ldots,n$.
	By Lemma \ref{lemma: T^k action on toric submanifolds}, there exists a fixed point of the $T^{k}$-action on $\overline{C_{\lambda}(V)}$.
\end{proof}
\end{prop}

We consider the case that $\{ \langle u^{\lambda}_{i},q_{j} \rangle \}_{i \not\in \mathcal{J}_{\lambda}}$ does not satisfy the condition $(1)$ in Lemma \ref{lemma: star conditions} for any $j=k+1,\ldots, n$, i.e., $\{ \langle u^{\lambda}_{i},q_{j} \rangle \}_{i \not\in \mathcal{J}_{\lambda}}$ satisfies the condition $(2)$ in Lemma \ref{lemma: star conditions} for any $j=k+1,\ldots, n$.

\begin{prop}
	\label{prop: clubsuit 2}
	Assume that the point $i_{V}^{\ast}(\lambda)$ is a vertex in the convex polytope $i_{V}^{\ast}(\Delta)$.
	If $\{ \langle u^{\lambda}_{i},q_{j} \rangle \}_{i \not\in \mathcal{J}_{\lambda}}$ satisfies the condition $(2)$ in Lemma \ref{lemma: star conditions} for any $j=k+1,\ldots, n$, then there exists a fixed point of the $T^{k}$-action on $\overline{C_{\lambda}(V)}$.
\begin{proof}
	From Lemma \ref{lemma: clubsuit 2}, for any $j$, there exist $i_{j}\neq i^{\prime}_{j} \not\in \mathcal{J}_{\lambda}$ such that 
	\begin{equation*}
		\{z^{\lambda} \in \varphi_{\lambda}(U_{\lambda})\mid z^{\lambda}_{i_{j}} = z^{\lambda}_{i^{\prime}_{j}} = 0, i_{j} \neq i^{\prime}_{j} \not\in \mathcal{J}_{\lambda} \} 
		\subset 
		\{z^{\lambda} \in \varphi_{\lambda}(U_{\lambda})\mid f^{\lambda}_{j}(z^{\lambda}) = 0 \} 
	\end{equation*}
	holds. Since for any $j$, 
	\begin{align*}
		&\{z^{\lambda} \in \varphi_{\lambda}(U_{\lambda}) \mid z^{\lambda}_{i} = 0, \; \text{for any} \; i\not\in \mathcal{J}_{\lambda}\} \\
		&\subset \{z^{\lambda} \in \varphi_{\lambda}(U_{\lambda})\mid z^{\lambda}_{i_{j}} = z^{\lambda}_{i^{\prime}_{j}} = 0,\; i_{j}\neq i^{\prime}_{j} \not\in \mathcal{J}_{\lambda} \}
	\end{align*}
	holds for some $i_{j} \neq i^{\prime}_{j} \not\in \mathcal{J}_{\lambda}$, we obtain 
	\begin{equation*}
		\{z^{\lambda} \in \varphi_{\lambda}(U_{\lambda}) \mid z^{\lambda}_{i} = 0, i\not\in \mathcal{J}_{\lambda}\}
		\subset \bigcap_{j=k+1}^{n}\{z^{\lambda} \in \varphi_{\lambda}(U_{\lambda}) \mid f^{\lambda}_{j}(z^{\lambda}) = 0 \}.
	\end{equation*}
	Since the right hand side is equal to $ \overline{C_{\lambda}(V)}$, we see that
	\begin{equation*}
		\{z^{\lambda} \in \varphi_{\lambda}(U_{\lambda}) \mid z^{\lambda}_{i} = 0, i\not\in \mathcal{J}_{\lambda}\}
	\subset \overline{C_{\lambda}(V)}.
	\end{equation*}
	In particular, as we noted that $\{z^{\lambda} \in \varphi_{\lambda}(U_{\lambda}) \mid z^{\lambda}_{i} = 0, \; i\not\in \mathcal{J}_{\lambda}\} \neq \emptyset$, we obtain 
	\begin{equation*}
		\{z^{\lambda} \in \varphi_{\lambda}(U_{\lambda}) \mid z^{\lambda}_{i} = 0, \; i\not\in \mathcal{J}_{\lambda}\} \cap \overline{C_{\lambda}(V)} \neq \emptyset.
	\end{equation*}
	By Lemma \ref{lemma: T^k action on toric submanifolds}, there exists a fixed point of the $T^{k}$-action on $\overline{C_{\lambda}(V)}$.
\end{proof}
\end{prop}

\begin{prop}
	\label{prop: interior iff no fixed points}
	If the point $i_{V}^{\ast}(\lambda)$ is a vertex in the convex polytope $i_{V}^{\ast}(\Delta)$, then there exists a fixed point of the $T^{k}$-action on $\overline{C_{\lambda}(V)}$.

\begin{proof}
	By Lemma \ref{lemma: star conditions}, if the point $i_{V}^{\ast}(\lambda)$ is a vertex in the convex polytope $i_{V}^{\ast}(\Delta)$, then $\{\langle u^{\lambda}_{i},q_{j} \rangle \}_{i \not\in \mathcal{J}_{\lambda}}$ satisfies either of the condition $(1)$ or the condition $(2)$ in Lemma \ref{lemma: star conditions} for each $j$.

	If there exists $j$ such that $\{\langle u^{\lambda}_{i},q_{j} \rangle \}_{i \not\in \mathcal{J}_{\lambda}}$ satisfies the condition $(1)$ in Lemma \ref{lemma: star conditions}, then by Proposition \ref{prop: clubsuit 1}, there exists a fixed point of the $T^{k}$-action on $\overline{C_{\lambda}(V)}$.
	If otherwise, i.e., if $\{\langle u^{\lambda}_{i},q_{j} \rangle \}_{i \not\in \mathcal{J}_{\lambda}}$ satisfies the condition $(2)$ in Lemma \ref{lemma: star conditions} for any $j$, then by Proposition \ref{prop: clubsuit 2}, there exists a fixed point of the $T^{k}$-action on $\overline{C_{\lambda}(V)}$.
\end{proof}
\end{prop}

By comparing the vertices in $\overline{\mu}(\overline{C(V)})$ with those in $i_{V}^{\ast}(\mu (X))$,
we say more about the image of the moment map $\overline{\mu}$.

\begin{thm}
\label{thm: moment polytope of toric submanifolds}
If $\overline{C(V)}$ is a complex submanifold in $X$, then we obtain 
$\overline{\mu}(\overline{C(V)}) = i_{V}^{\ast}(\mu (X))$ in $(\mathfrak{t}^{k})^{\ast}$.

\begin{proof}
	Since the map $\overline{\mu}$ is the moment map, the image of $\overline{\mu}$ is the convex hull of the images of the fixed points of the $T^{k}$-action on $\overline{C(V)}$.
	We classified the fixed points of the $T^{k}$-action on $X$ (Lemma \ref{lemma: fixed points of subtorus action}) and those of the $T^{k}$-action on $\overline{C(V)}$ (Lemma \ref{lemma: T^k action on toric submanifolds}).

	Since $\overline{C(V)} \subset X$, we obtain $\overline{\mu}(\overline{C(V)}) \subset i_{V}^{\ast}(\mu(X)) = i_{V}^{\ast}(\Delta)$.
	In particular, by Lemma \ref{lemma: fixed point in face}, if $z^{\lambda} \in \overline{C_{\lambda}(V)}$ is a fixed point of the $T^{k}$-action on $\overline{C_{\lambda}(V)}$, then $\overline{\mu}(z^{\lambda}) = i_{V}^{\ast}(\lambda) \in i_{V}^{\ast}(\Delta)$ for the vertex $\lambda$.

	Since Proposition \ref{prop: interior iff no fixed points} shows that if $i_{V}^{\ast}(\lambda)$ is a vertex of $i_{V}^{\ast}(\Delta)$, then there exists a fixed point $z^{\lambda}$ of the $T^{k}$-action on $\overline{C_{\lambda}(V)}$ such that $\overline{\mu}(z^{\lambda})= i_{V}^{\ast}(\lambda)$.

	Thus, the set of the vertices of $\overline{\mu}(\overline{C(V)})$ coincides with the set of the vertices of $i_{V}^{\ast}(\Delta)$.
	Since the map $\overline{\mu}$ is a moment map for the $T^{k}$-action on $\overline{C(V)}$, by the convexity theorem \cite{MR642416,MR664117}, the image of $\overline{\mu}$ is the convex hull of the images of the fixed points of the $T^{k}$-action on $\overline{C(V)}$.
	Since $i_{V}^{\ast}(\Delta)$ is the convex hull of the images of the vertices of $\Delta$ by $i_{V}^{\ast}$, we obtain $\overline{\mu}(\overline{C(V)}) = i_{V}^{\ast}(\Delta)$.
\end{proof}
\end{thm}

We say a submanifold $\overline{C(V)}$ to be a \textit{torus-equivariantly embedded toric manifold} in a toric manifold $X$.

\section{Examples of Torus-equivariantly Embedded Toric Manifolds}
\label{sec: example}

We demonstrate examples of $\overline{C(V)}$ and check whether they are torus-equivariantly embedded toric manifolds or not. When $\overline{C(V)}$ is smooth, we further draw figures of $\overline{D(V)}:= \mu (\overline{C(V)})$ for each example.

\subsection{Examples of Torus-equivariantly Embedded Toric Manifolds in $\mathbb{C}P^{2}$}
\label{subsec: example}

We give examples for $\overline{C(V)}$ and check whether $\overline{C(V)}$ is a complex submanifold in $X = \mathbb{C}P^{2}$.

Delzant polytopes of $\mathbb{C}P^{2}$ are isosceles right triangles. 
As in Example \ref{eg: cp2, o.k.} and Example \ref{eg: cp2, n.g.},
define the points $\lambda$, $\mu$, $\sigma$ in $(\mathfrak{t}^{2})^{\ast} \cong \mathbb{R}^{2}$ by 
	\begin{equation*}
		\lambda = (0,0), ~
		\mu = (2,0), ~ 
		\sigma = (0,2).
	\end{equation*}
Let $\Delta$ be a polytope defined by the convex hull of the points $\lambda$, $\mu$, $\sigma$.
We define the inward pointing normal vectors to the facets by 
\begin{equation*}
	u^{\lambda}_{1} =
	\left[
		\begin{matrix}
			1 \\ 0
		\end{matrix}
	\right],
	u^{\lambda}_{2} =
	\left[
		\begin{matrix}
			0 \\ 1
		\end{matrix}
	\right], 
	u^{\mu}_{1} = 
	\left[
		\begin{matrix}
			-1 \\ -1
		\end{matrix}
	\right],
	u^{\mu}_{2} = 
	\left[
		\begin{matrix}
			1 \\ 0
		\end{matrix}
	\right],
	u^{\sigma}_{1} = 
	\left[
		\begin{matrix}
			0 \\ 1
		\end{matrix}
	\right],
	u^{\sigma}_{2} = 
	\left[
		\begin{matrix}
			-1 \\ -1
		\end{matrix}
	\right],
\end{equation*}
where $\Lambda = \{ \lambda,\mu, \sigma\}$.

\begin{figure}[hbtp]
	\centering
	\begin{tabular}{cc}
	\begin{minipage}[t]{0.5\hsize}
		\centering
		\includegraphics[keepaspectratio,width=5cm]{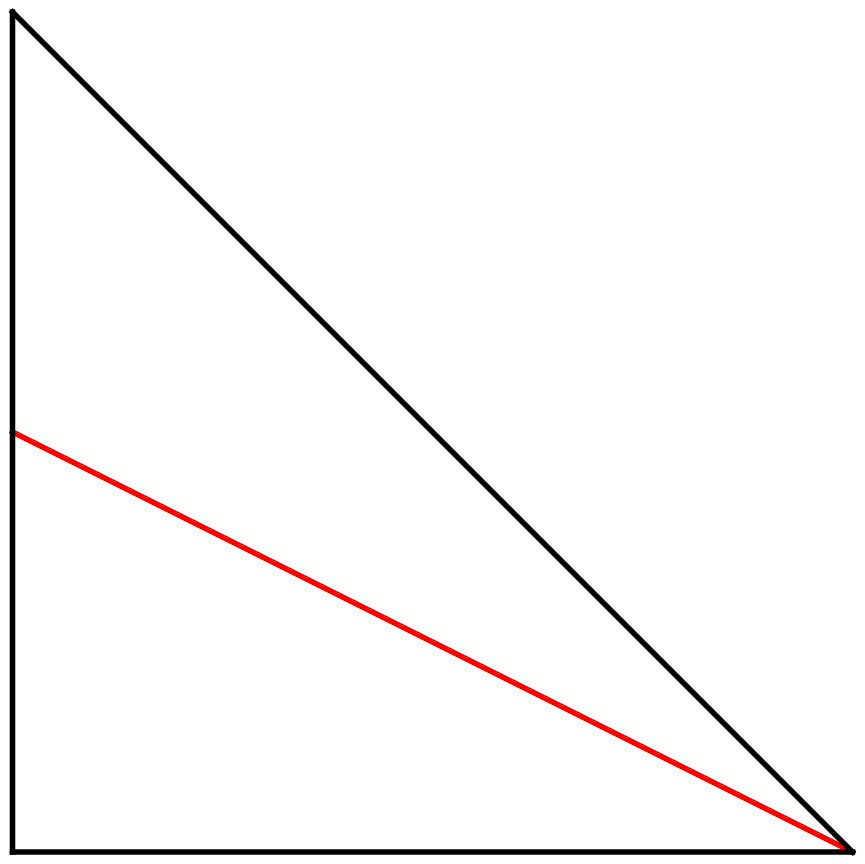}
		\caption{$\overline{D(V)}$ in Example \ref{eg: example10cp2} when $a = (0,0)$}
		\label{fig: example10cp2}
	\end{minipage} &
	\begin{minipage}[t]{0.5\hsize}
		\centering
		\includegraphics[keepaspectratio,width=5cm]{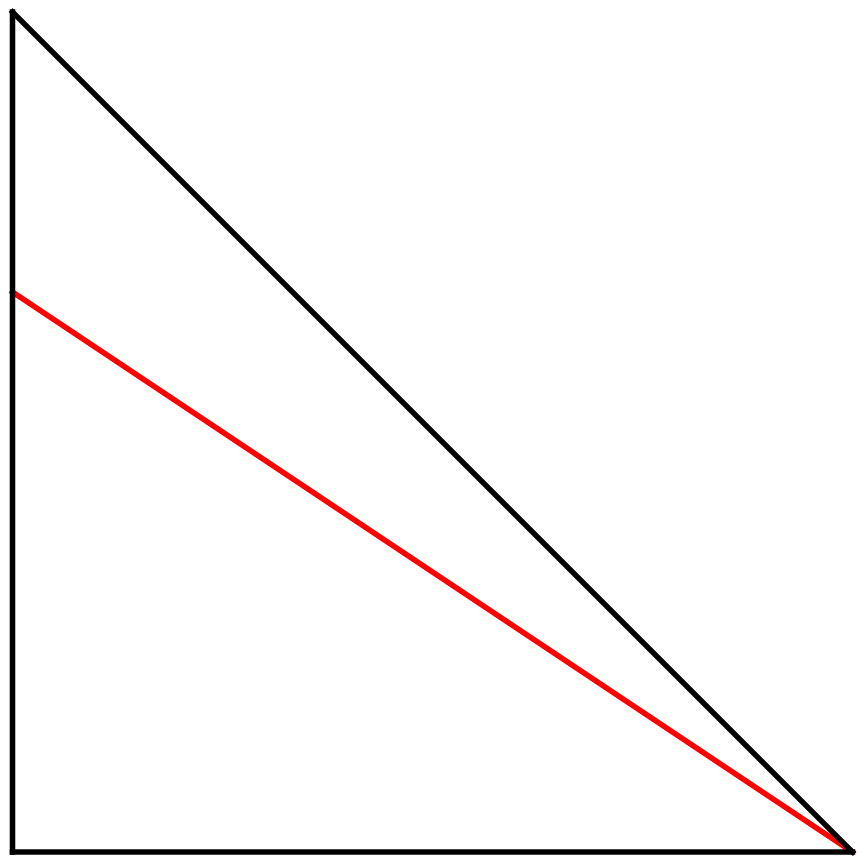}
		\caption{$\overline{D(V)}$ in Example \ref{eg: example10cp2} when $a = (0,\log 2)$}
		\label{fig: example10acp2}
	\end{minipage}
\end{tabular}
\end{figure}

\begin{example}
	\label{eg: example10cp2}
	Let $X = \mathbb{C}P^{2}$, $k = 1$, and $V = \mathbb{R}p + a \; (a \in \mathbb{R}^{2})$ be an affine subspace spanned by $p = {^t[1 ~ 0]}$.
	Then, we can choose a basis $q$ of the orthogonal subspace to the linear part of $V$ as $q = {^t [0 ~ 1]}$. 
	In this case, $\overline{C(V)}$ is a complex submanifold in $X$. 
	Indeed, we give $f^{\lambda},f^{\mu},f^{\sigma}$ by 
	\begin{equation*}
		f^{\lambda} = z^{\lambda}_{2}-e^{\langle a, e_{2} \rangle}, \;
		f^{\mu} = 1- e^{\langle a, e_{2} \rangle}z^{\mu}_{1}, \;
		f^{\sigma} = z^{\sigma}_{1} - e^{\langle a, e_{2} \rangle}z^{\sigma}_{2},
	\end{equation*}
	respectively. Since the Jacobian matrices are expressed as
	\begin{equation*}
		Df^{\lambda} = 
		[\begin{matrix}
			0 & 1
		\end{matrix}], \;
		Df^{\mu} = 
		[\begin{matrix}
			-e^{\langle a, e_{2} \rangle} & 0
		\end{matrix}], \;
		Df^{\sigma} = 
		[\begin{matrix}
			1 & -e^{\langle a, e_{2} \rangle}
		\end{matrix}],
	\end{equation*}
	respectively, we see the rank of each matrix is one.

	Figure \ref{fig: example10cp2} describes $\overline{D(V)}$ when $a =(0,0)$.
	Figure \ref{fig: example10acp2} describes $\overline{D(V)}$ when $a = (0,\log 2)$.
\end{example}

\begin{figure}[hbtp]
	\centering
	\begin{tabular}{cc}
	\begin{minipage}[t]{0.5\hsize}
		\centering
		\includegraphics[keepaspectratio,width=5cm]{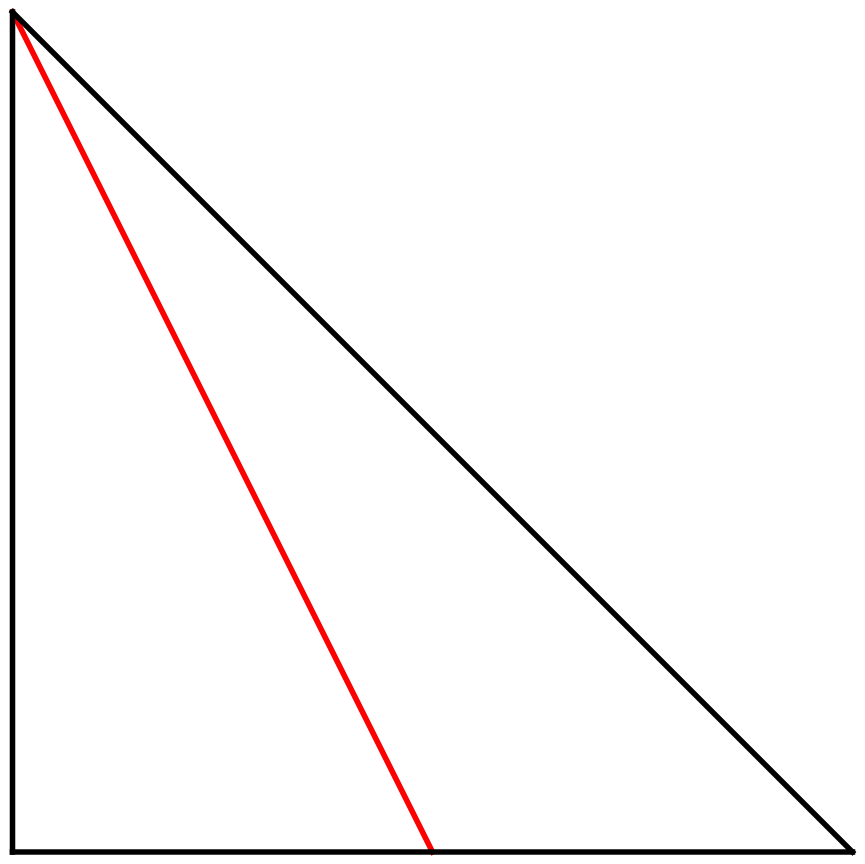}
		\caption{$\overline{D(V)}$ in Example \ref{eg: example01cp2} when $a = (0,0)$}
		\label{fig: example01cp2}
	\end{minipage} &
	\begin{minipage}[t]{0.5\hsize}
		\centering
		\includegraphics[keepaspectratio,width=5cm]{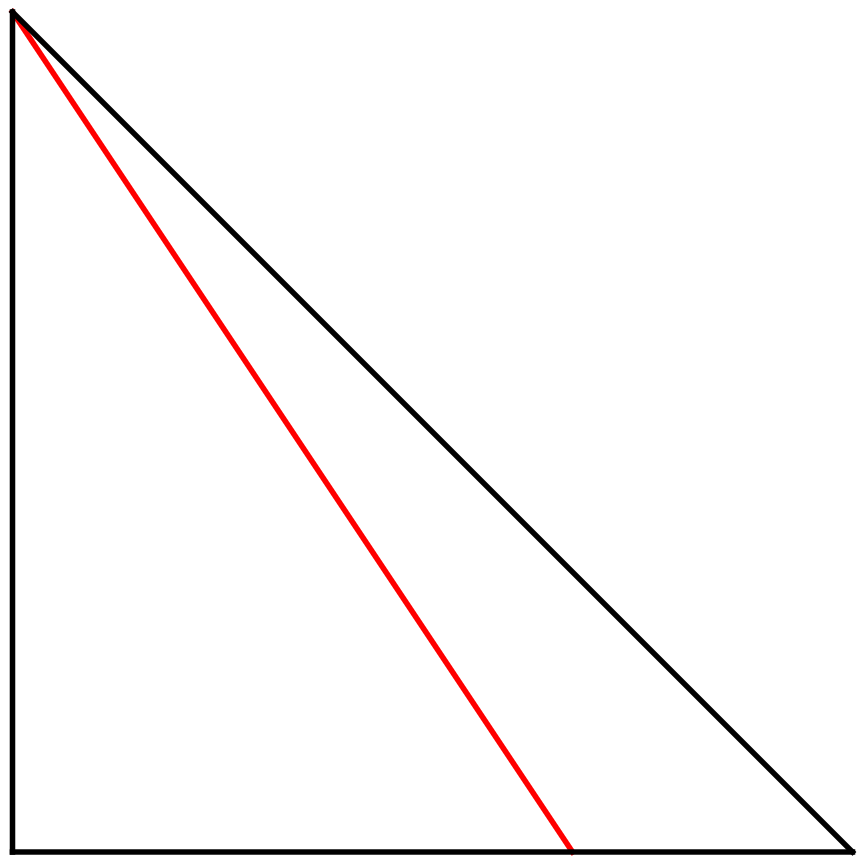}
		\caption{$\overline{D(V)}$ in Example \ref{eg: example01cp2} when $a = (\log 2,0)$}
		\label{fig: example01acp2}
	\end{minipage}
\end{tabular}
\end{figure}

\begin{example}
	\label{eg: example01cp2}
	Let $X = \mathbb{C}P^{2}$, $k = 1$, and $V = \mathbb{R}p + a \; (a \in \mathbb{R}^{2})$ be an affine subspace spanned by $p = {^t[0 ~ 1]}$.
	Then, we can choose a basis $q$ of the orthogonal subspace to the linear part of $V$ as $q = {^t [1 ~ 0]}$. 
	In this case, $\overline{C(V)}$ is a complex submanifold in $X$. 
	Indeed, we give $f^{\lambda},f^{\mu},f^{\sigma}$ by 
	\begin{equation*}
		f^{\lambda} = z^{\lambda}_{1}-e^{\langle a,e_{1}\rangle}, \;
		f^{\mu} = z^{\mu}_{2} - e^{\langle a,e_{1}\rangle}z^{\mu}_{1}, \;
		f^{\sigma} = 1 - e^{\langle a,e_{1}\rangle}z^{\sigma}_{2},
	\end{equation*}
	respectively. Since the Jacobian matrices are expressed as
	\begin{equation*}
		Df^{\lambda} = [
			\begin{matrix}
				1 & 0
			\end{matrix}
		], \;
		Df^{\mu} = [
			\begin{matrix}
				-e^{\langle a,e_{1}\rangle} & 1
			\end{matrix}
		], \;
		Df^{\sigma} = [
			\begin{matrix}
				0 & -e^{\langle a,e_{1}\rangle}
			\end{matrix}
		],
	\end{equation*}
	respectively, we see the rank of each matrix is one.
	
	Figure \ref{fig: example01cp2} describes $\overline{D(V)}$ when $a =(0,0)$.
	Figure \ref{fig: example01acp2} describes $\overline{D(V)}$ when $a = (\log 2,0)$.
\end{example}

\begin{figure}[hbtp]
	\centering
	\begin{tabular}{cc}
	\begin{minipage}[t]{0.5\hsize}
		\centering
		\includegraphics[keepaspectratio,width=5cm]{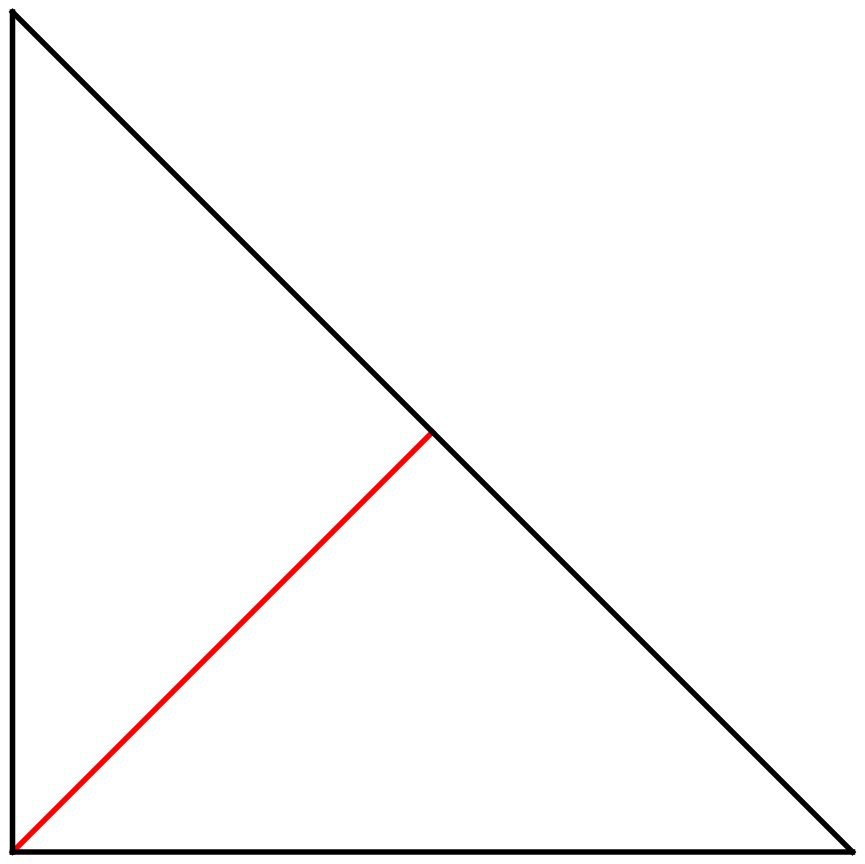}
		\caption{$\overline{D(V)}$ in Example \ref{eg: example11cp2} when $a= (0,0)$}
		\label{fig: example11cp2}
	\end{minipage} &
	\begin{minipage}[t]{0.5\hsize}
		\centering
		\includegraphics[keepaspectratio,width=5cm]{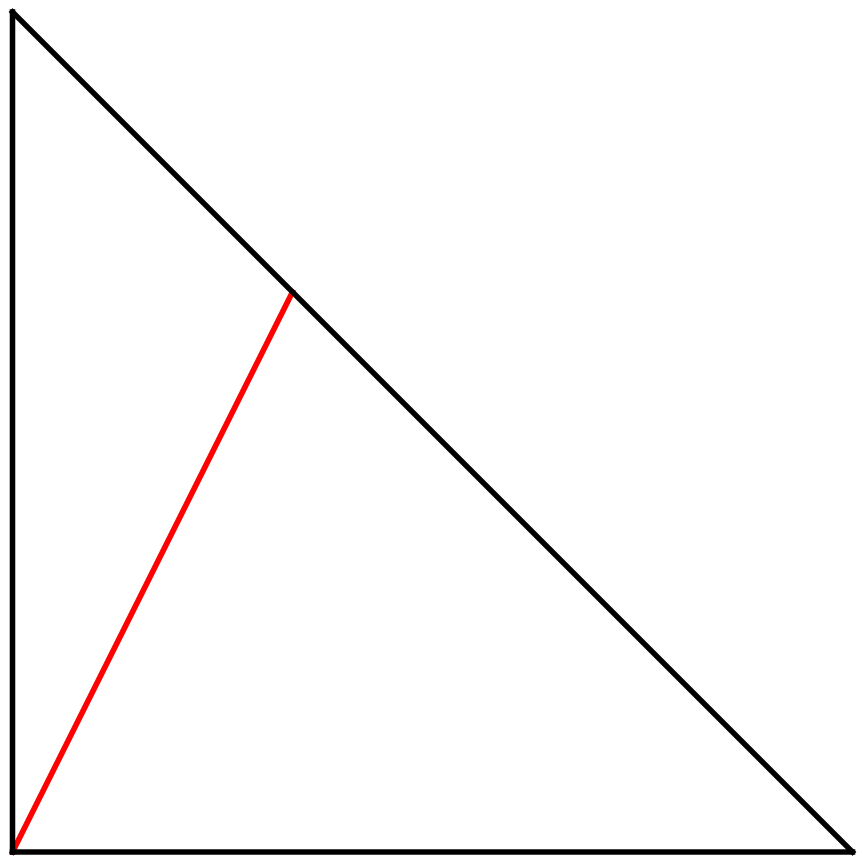}
		\caption{$\overline{D(V)}$ in Example \ref{eg: example11cp2} when $a = (0,\log 2)$}
		\label{fig: example11acp2}
	\end{minipage}
	\end{tabular}
\end{figure}

\begin{example}
	\label{eg: example11cp2}
	Let $X = \mathbb{C}P^{2}$, $k = 1$, and $V = \mathbb{R}p + a \; (a \in \mathbb{R}^{2})$ be an affine subspace spanned by $p = {^t[1 ~ 1]}$.
	Then, we can choose a basis $q$ of the orthogonal subspace to the linear part of $V$ as $q = {^t [1 ~ -1]}$. 
	In this case, we show in Example \ref{eg: cp2, o.k.} that $\overline{C(V)}$ is a complex submanifold in $X$ for $a = 0$. By similar calculation, we see that $\overline{C(V)}$ is a complex submanifold in $X$ for arbitrary $a \in \mathbb{R}^{2}$.

	Figure \ref{fig: example11cp2} describes $\overline{D(V)}$ when $a=(0,0)$.
	Figure \ref{fig: example11acp2} describes $\overline{D(V)}$ when $a = (0,\log 2)$.
\end{example}

\begin{figure}[hbtp]
	\centering
	\begin{tabular}{cc}
	\begin{minipage}[t]{0.5\hsize}
		\centering
		\includegraphics[keepaspectratio,width=5cm]{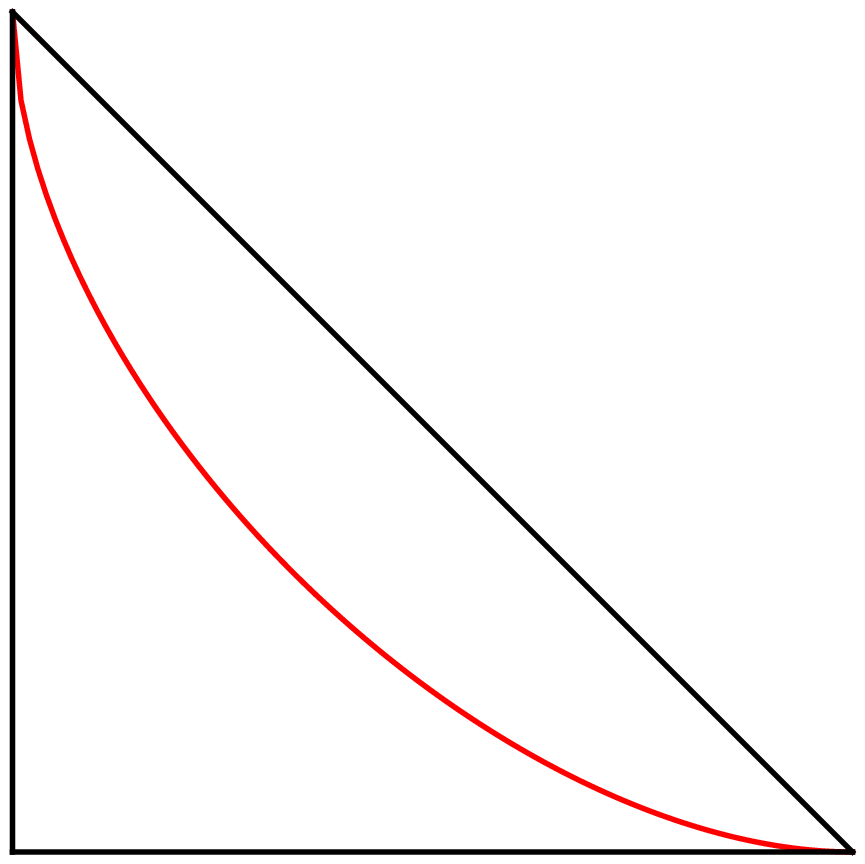}
		\caption{$\overline{D(V)}$ in Example \ref{eg: example1-1cp2} when $a = (0,0)$}
		\label{fig: example1-1cp2}
	\end{minipage} &
	\begin{minipage}[t]{0.5\hsize}
		\centering
		\includegraphics[keepaspectratio,width=5cm]{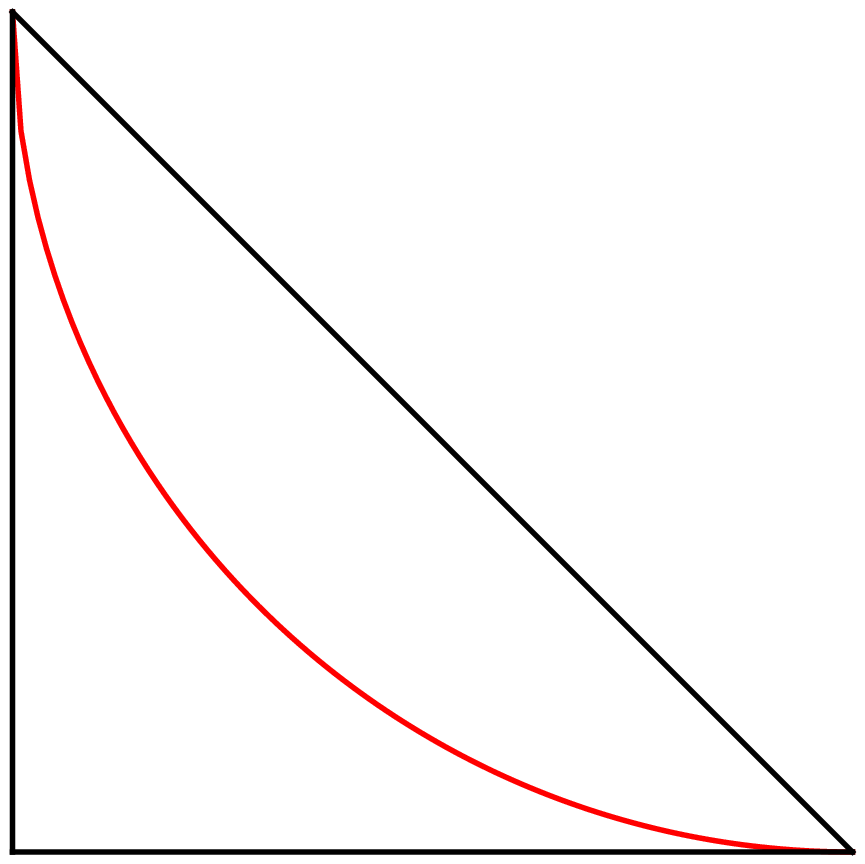}
		\caption{$\overline{D(V)}$ in Example \ref{eg: example1-1cp2} when $a = (-\log 2,0)$}
		\label{fig: example1-1acp2}
	\end{minipage}
	\end{tabular}
\end{figure}

\begin{example}
	\label{eg: example1-1cp2}
	Let $X = \mathbb{C}P^{2}$, $k = 1$, and $V = \mathbb{R}p + a \; (a \in \mathbb{R}^{2})$ be an affine subspace spanned by $p = {^t[1 ~ -1]}$.
	Then, we can choose a basis $q$ of the orthogonal subspace to the linear part of $V$ as $q = {^t [1 ~ 1]}$. 
	In this case, $\overline{C(V)}$ is a complex submanifold in $X$. 
	Indeed, we give $f^{\lambda},f^{\mu},f^{\sigma}$ by 
	\begin{align*}
		f^{\lambda} &= z^{\lambda}_{1}z^{\lambda}_{2} -e^{\langle a,e_{1} \rangle + \langle a,e_{2} \rangle}, \\
		f^{\mu} &= z^{\mu}_{2} - e^{\langle a,e_{1} \rangle + \langle a,e_{2} \rangle}(z^{\mu}_{1})^{2}, \\
		f^{\sigma} &= e^{\langle a,e_{1} \rangle + \langle a,e_{2} \rangle}(z^{\sigma} - (z^{\sigma}_{2})^{2}),
	\end{align*}
	respectively. 
	Note that $(0,0) \notin \overline{C_{\lambda}(V)}$.
	Since the Jacobian matrices are expressed as
	\begin{align*}
		Df^{\lambda} &=[
			\begin{matrix}
				z^{\lambda}_{2} & z^{\lambda}_{1}
			\end{matrix}
		], \\
		Df^{\mu} &= [
			\begin{matrix}
				-2 e^{\langle a,e_{1} \rangle + \langle a,e_{2} \rangle} z^{\mu}_{1} & 1
			\end{matrix}
		], \\
		Df^{\sigma} &= [
			\begin{matrix}
				1 & -2 e^{\langle a,e_{1} \rangle + \langle a,e_{2} \rangle} z^{\sigma}_{2}
			\end{matrix}
		],
	\end{align*}
	respectively, we see the rank of each matrix is one.
	
	Figure \ref{fig: example1-1cp2} describes $\overline{D(V)}$ when $a = (0,0)$.
	Figure \ref{fig: example1-1acp2} describes $\overline{D(V)}$ when $a = (-\log 2,0)$.
\end{example}

\begin{figure}[hbtp]
	\centering
	\begin{tabular}{cc}
	\begin{minipage}[t]{0.5\hsize}
		\centering
		\includegraphics[keepaspectratio,width=5cm]{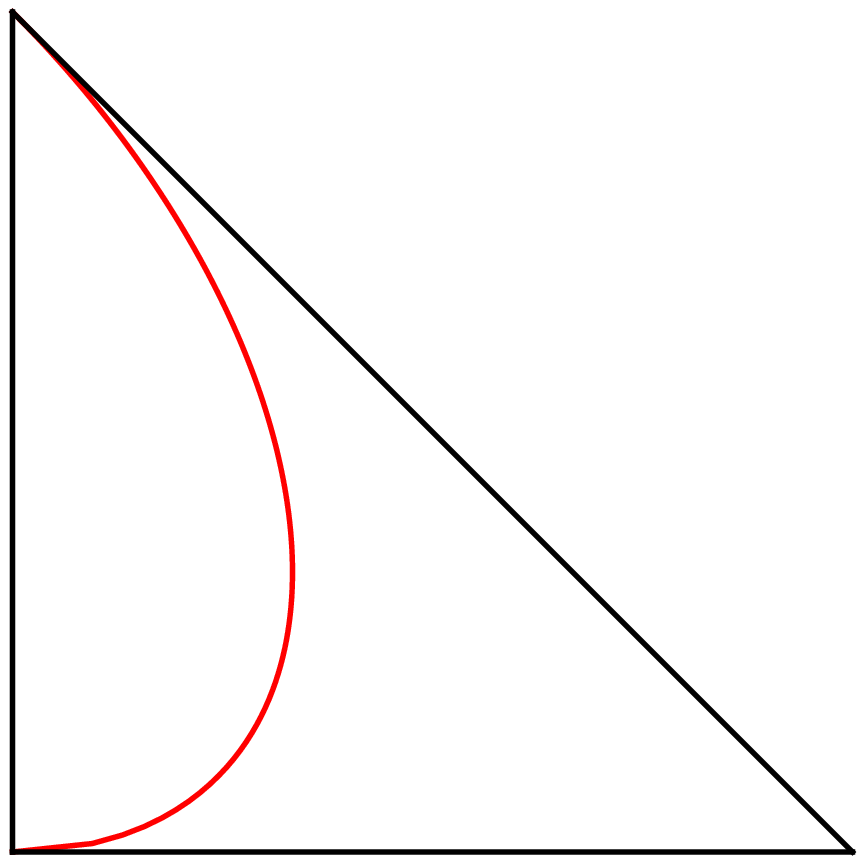}
		\caption{$\overline{D(V)}$ in Example \ref{eg: example12cp2} when $a = (0, 0)$}
		\label{fig: example12cp2}
	\end{minipage} &
	\begin{minipage}[t]{0.5\hsize}
		\centering
		\includegraphics[keepaspectratio,width=5cm]{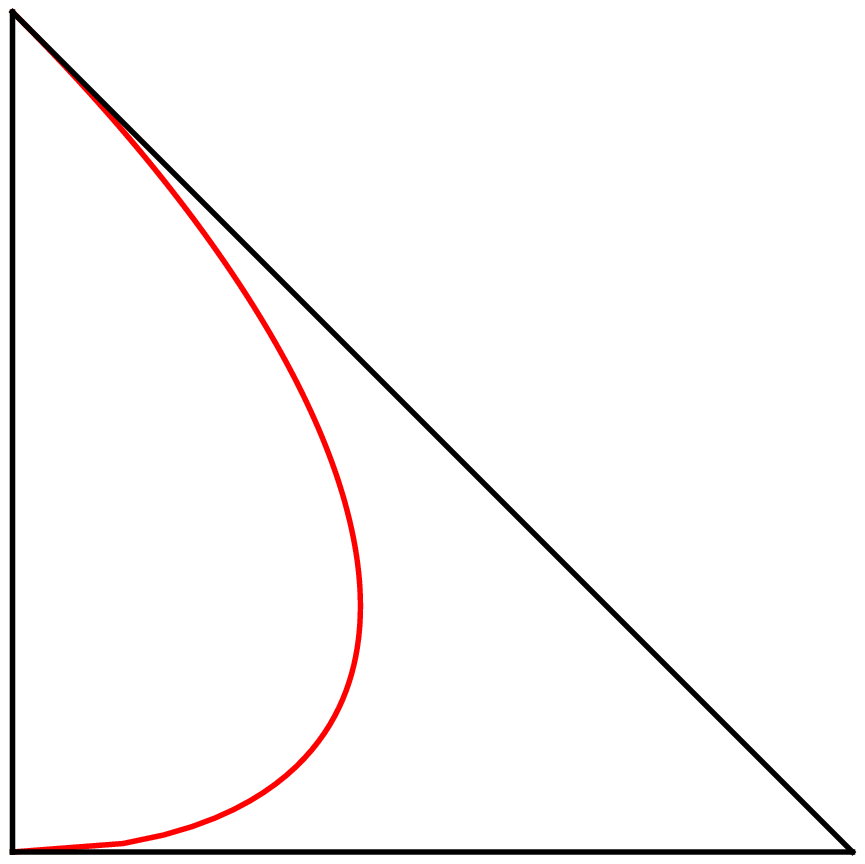}
		\caption{$\overline{D(V)}$ in Example \ref{eg: example12cp2} when $a = (0, -\log 2)$}
		\label{fig: example12acp2}
	\end{minipage}
	\end{tabular}
\end{figure}

\begin{example}
	\label{eg: example12cp2}
	Let $X = \mathbb{C}P^{2}$, $k = 1$, and $V = \mathbb{R}p + a \; (a \in \mathbb{R}^{2})$ be an affine subspace spanned by $p = {^t[1 ~ 2]}$.
	Then, we can choose a basis $q$ of the orthogonal subspace to the linear part of $V$ as $q = {^t [2 ~ -1]}$. 
	In this case, $\overline{C(V)}$ is a complex submanifold in $X$. 
	Indeed, we give $f^{\lambda},f^{\mu},f^{\sigma}$ by 
	\begin{align*}
		f^{\lambda} 
		&= 
		(z^{\lambda}_{1})^{2} - e^{2\langle a,e_{1} \rangle - \langle a,e_{2} \rangle} z^{\lambda}_{2}, \\ 
		f^{\mu} 
		&= 
		(z^{\mu}_{2})^{2} - e^{2\langle a,e_{1} \rangle - \langle a,e_{2} \rangle} z^{\mu}_{1},\\
		f^{\sigma} 
		&= 
		1 - e^{2\langle a,e_{1} \rangle - \langle a,e_{2} \rangle}z^{\sigma}_{1}z^{\sigma}_{2},
	\end{align*}
	respectively. 
	Note that $(0,0) \notin \overline{C_{\sigma}(V)}$. 
	Since the Jacobian matrices are expressed as
	\begin{align*}
		Df^{\lambda} &= [
			\begin{matrix}
				2 z^{\lambda} & - e^{2\langle a,e_{1} \rangle - \langle a,e_{2} \rangle}
			\end{matrix}
		],\\
		Df^{\mu} &= [
			\begin{matrix}
				- e^{2\langle a,e_{1} \rangle - \langle a,e_{2} \rangle} & 2 z^{\mu}_{2}
			\end{matrix}
		], \\
		Df^{\sigma} &= 
		-e^{2\langle a,e_{1} \rangle - \langle a,e_{2} \rangle}
		[
			\begin{matrix}
				z^{\sigma}_{2} & z^{\sigma}_{1}
			\end{matrix}
		],
	\end{align*}
	respectively, we see the rank of each matrix is one.
	
	Figure \ref{fig: example12cp2} describes $\overline{D(V)}$ when $a = (0,0)$.
	Figure \ref{fig: example12acp2} describes $\overline{D(V)}$ when $a = (0, -\log 2)$.
\end{example}

\begin{figure}[hbtp]
	\centering
	\begin{tabular}{cc}
	\begin{minipage}[t]{0.5\hsize}
		\centering
		\includegraphics[keepaspectratio,width=5cm]{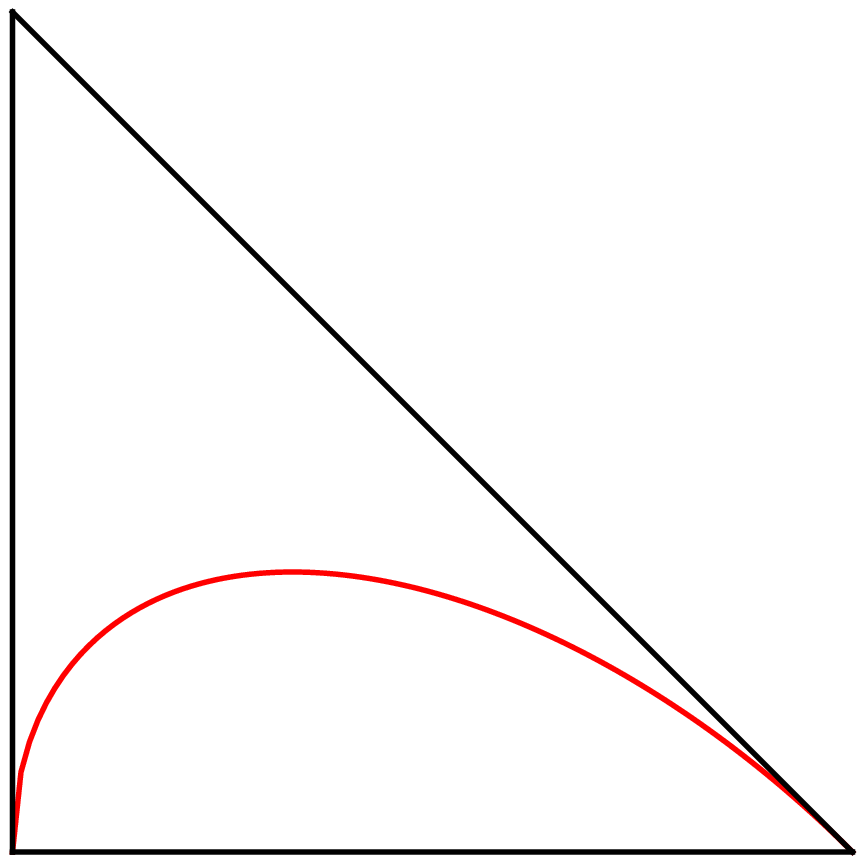}
		\caption{$\overline{D(V)}$ in Example \ref{eg: example21cp2} when $a=(0,0)$}
		\label{fig: example21cp2}
	\end{minipage} &
	\begin{minipage}[t]{0.5\hsize}
		\centering
		\includegraphics[keepaspectratio,width=5cm]{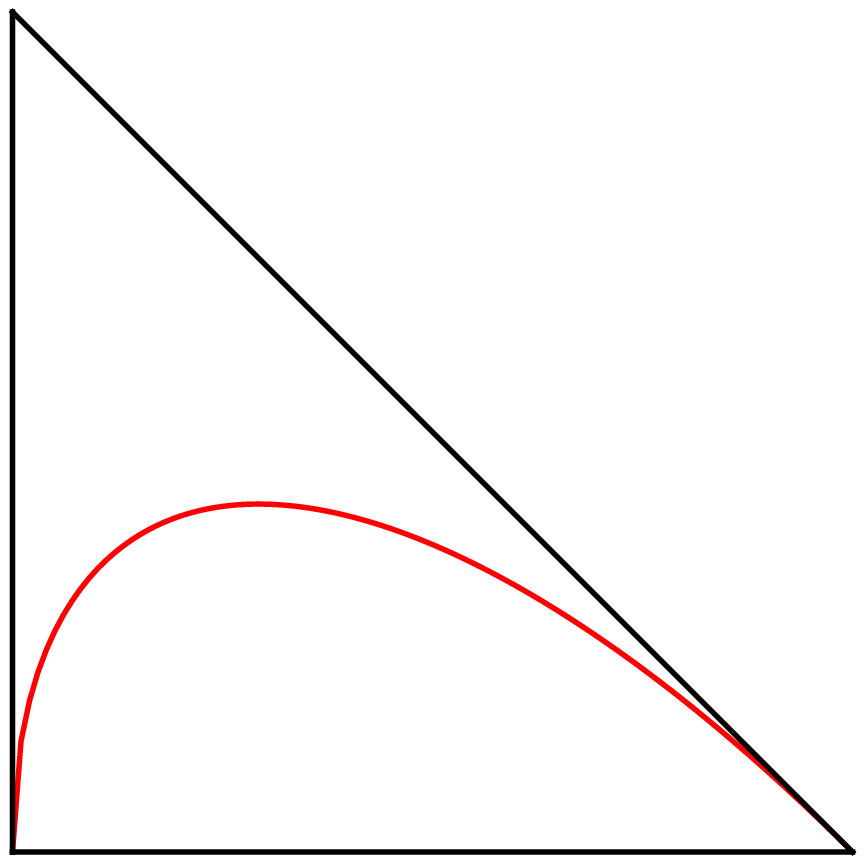}
		\caption{$\overline{D(V)}$ in Example \ref{eg: example21cp2} when $a=(-\log 2,0)$}
		\label{fig: example21acp2}
	\end{minipage}
	\end{tabular}
\end{figure}

\begin{example}
	\label{eg: example21cp2}
	Let $X = \mathbb{C}P^{2}$, $k = 1$, and $V = \mathbb{R}p + a ~(a \in \mathbb{R}^{2})$ be an affine subspace spanned by $p = {^t[2 ~ 1]}$.
	Then, we can choose a basis $q$ of the orthogonal subspace to the linear part of $V$ as $q = {^t [1 ~ -2]}$. 
	In this case, $\overline{C(V)}$ is a complex submanifold in $X$. 
	Indeed, we give $f^{\lambda},f^{\mu},f^{\sigma}$ by 
	\begin{align*}
		f^{\lambda} 
		&= 
		z^{\lambda}_{1} - e^{\langle a,e_{1} \rangle - 2 \langle a,e_{2} \rangle}(z^{\lambda}_{2})^{2}, \\
		f^{\mu} 
		&= 
		z^{\mu}_{1} z^{\mu}_{2} - e^{\langle a,e_{1} \rangle - 2 \langle a,e_{2} \rangle}, \\
		f^{\sigma} 
		&=
		z^{\sigma}_{2} -e^{\langle a,e_{1} \rangle - 2 \langle a,e_{2} \rangle}(z^{\sigma}_{1})^{2},
	\end{align*}
	respectively. 
	Note that $(0,0) \notin \overline{C_{\mu}(V)}$. 
	Since the Jacobian matrices are expressed as
	\begin{align*}
		Df^{\lambda} &= [
			\begin{matrix}
				1 & - 2e^{\langle a,e_{1} \rangle - 2 \langle a,e_{2} \rangle}z^{\lambda}_{2}
			\end{matrix}
		], \\
		Df^{\mu} &= [
			\begin{matrix}
				z^{\mu}_{2} & z^{\mu}_{1}
			\end{matrix}
		], \\
		Df^{\sigma} &= [
			\begin{matrix}
				-2 e^{\langle a,e_{1} \rangle - 2 \langle a,e_{2} \rangle} z^{\sigma}_{1} & 1
			\end{matrix}
		],
	\end{align*}
	respectively, we see the rank of each matrix is one.
	
	Figure \ref{fig: example21cp2} describes $\overline{D(V)}$ when $a = (0,0)$.
	Figure \ref{fig: example21acp2} describes $\overline{D(V)}$ when $a = (-\log 2, 0)$.
\end{example}

In the following examples, we treat $\overline{C(V)}$ which does not become a complex submanifold in $\mathbb{C}P^{2}$. 

\begin{example}
	\label{eg: nonexample1Acp2}
	Let $X = \mathbb{C}P^{2}$, $k = 1$, and $V = \mathbb{R}p + a \; (a \in \mathbb{R}^{2})$ be an affine subspace spanned by $p = {^t[1 ~ \alpha]}$ for all integers $\alpha$ greater than or equal to three.
	Then, we can choose a basis $q$ of the orthogonal subspace to $V$ as $q = {^t [\alpha ~ -1]}$. 
	In this case, $\overline{C(V)}$ is not a complex submanifold in $X$. 
	Indeed, we give $f^{\lambda},f^{\mu},f^{\sigma}$ by 
	\begin{align*}
		f^{\lambda} 
		&= 
		(z^{\lambda}_{1})^{\alpha} - e^{\alpha \langle a,e_{1} \rangle - \langle a,e_{2} \rangle} z^{\lambda}_{2}, \\
		f^{\mu} 
		&= 
		(z^{\mu}_{2})^{\alpha} - e^{\alpha \langle a,e_{1} \rangle - \langle a,e_{2} \rangle} (z^{\mu}_{1})^{\alpha -1}, \\
		f^{\sigma} 
		&= 
		1 - e^{\alpha \langle a,e_{1} \rangle - \langle a,e_{2} \rangle} z^{\sigma}_{1} (z^{\sigma}_{2})^{\alpha -1},
	\end{align*}
	respectively. 
	Note that $(0,0) \notin \overline{C_{\sigma}(V)}$.
	Since the Jacobian matrices are expressed as
	\begin{align*}
		Df^{\lambda} &= [
			\begin{matrix}
				\alpha (z^{\lambda}_{1})^{\alpha - 1} & - e^{\alpha \langle a,e_{1} \rangle - \langle a,e_{2} \rangle}
			\end{matrix}
		], \\
		Df^{\mu} &= [
			\begin{matrix}
				- (\alpha - 1)e^{\alpha \langle a,e_{1} \rangle - \langle a,e_{2} \rangle} (z^{\mu}_{1})^{\alpha - 2} & \alpha (z^{\mu}_{2})^{\alpha -1}
			\end{matrix}
		], \\
		Df^{\sigma} &= 
		-e^{\alpha \langle a,e_{1} \rangle - \langle a,e_{2} \rangle}
		[
			\begin{matrix}
				(z^{\sigma}_{2})^{\alpha -1} & (\alpha -1)z^{\sigma}_{1}(z^{\sigma}_{2})^{\alpha -2}
			\end{matrix}
		],
	\end{align*}
	respectively, we see the rank of $Df^{\lambda}$ and $Df^{\sigma}$ is one. However, the rank of $Df^{\mu}$ becomes zero when $(z^{\mu}_{1},z^{\mu}_{2}) = (0,0) \in \overline{C_{\mu}(V)}$. 
\end{example}

\begin{example}
	\label{eg: nonexampleA1cp2}
	Let $X = \mathbb{C}P^{2}$, $k = 1$, and $V = \mathbb{R}p + a \; (a \in \mathbb{R}^{2})$ be an affine subspace spanned by $p = {^t[\alpha ~ 1]}$ for all integers $\alpha$ greater than or equal to three.
	Then, we can choose a basis $q$ of the orthogonal subspace to $V$ as $q = {^t [1 ~ -\alpha]}$. 
	In this case, $\overline{C(V)}$ is not a complex submanifold in $X$. 
	Indeed, we give $f^{\lambda},f^{\mu},f^{\sigma}$ by 
	\begin{align*}
		f^{\lambda} 
		&= 
		z^{\lambda}_{1} - e^{\langle a,e_{1} \rangle - \alpha \langle a,e_{2} \rangle}(z^{\lambda}_{2})^{\alpha}, \\
		f^{\mu} 
		&= 
		(z^{\mu}_{1})^{\alpha -1}z^{\mu}_{2} - e^{\langle a,e_{1} \rangle - \alpha \langle a,e_{2} \rangle}, \\
		f^{\sigma} 
		&= 
		(z^{\sigma}_{2})^{\alpha -1} -e^{\langle a,e_{1} \rangle - \alpha \langle a,e_{2} \rangle}(z^{\sigma}_{1})^{\alpha},
	\end{align*}
	respectively. 
	Note that $(0,0) \notin \overline{C_{\mu}(V)}$.
	Since the Jacobian matrices are expressed as
	\begin{align*}
		Df^{\lambda} 
		&= [
			\begin{matrix}
				1 & - \alpha e^{\langle a,e_{1} \rangle - \alpha \langle a,e_{2} \rangle} (z^{\lambda}_{2})^{\alpha -1}
			\end{matrix}
		], \\
		Df^{\mu} 
		&= [
			\begin{matrix}
				(\alpha -1)(z^{\mu}_{1})^{\alpha -2}z^{\mu}_{2} & (z^{\mu}_{1})^{\alpha -1}
			\end{matrix}	
		], \\
		Df^{\sigma} 
		&= [
			\begin{matrix}
				-\alpha e^{\langle a,e_{1} \rangle - \alpha \langle a,e_{2} \rangle} (z^{\sigma}_{1})^{\alpha -1} & (\alpha -1) (z^{\sigma}_{2})^{\alpha -2}
			\end{matrix}
		],
	\end{align*}
	respectively, we see the rank of $Df^{\lambda}$ and $Df^{\sigma}$ is one. However, the rank of $Df^{\sigma}$ becomes zero at the point $(z^{\sigma}_{1},z^{\sigma}_{2}) = (0,0) \in \overline{C_{\sigma}(V)}$. 
\end{example}

\begin{example}
	\label{eg: nonexample1-Acp2}
	Let $X = \mathbb{C}P^{2}$, $k = 1$, and $V = \mathbb{R}p + a \; (a \in \mathbb{R}^{2})$ be an affine subspace spanned by $p = {^t[1 ~ -\alpha]}$ for all integers $\alpha$ greater than or equal to two.
	Then, we can choose a basis $q$ of the orthogonal subspace to $V$ as $q = {^t [\alpha ~ 1]}$. 
	In this case, $\overline{C(V)}$ is not a complex submanifold in $X$. 
	Indeed, we give $f^{\lambda},f^{\mu},f^{\sigma}$ by 
	\begin{align*}
		f^{\lambda} 
		&= 
		(z^{\lambda}_{1})^{\alpha}z^{\lambda}_{2} - e^{\alpha \langle a,e_{1} \rangle + \langle a,e_{2} \rangle}, \\
		f^{\mu} 
		&= 
		(z^{\mu}_{2})^{\alpha} - e^{\alpha \langle a,e_{1} \rangle + \langle a,e_{2} \rangle} (z^{\mu}_{1})^{\alpha + 1}, \\
		f^{\sigma} 
		&= 
		z^{\sigma}_{1} - e^{\alpha \langle a,e_{1} \rangle + \langle a,e_{2} \rangle} (z^{\sigma}_{2})^{\alpha +1},
	\end{align*}
	respectively. 
	Note that $(0,0) \notin \overline{C_{\lambda}(V)}$.
	Since the Jacobian matrices are expressed as
	\begin{align*}
		Df^{\lambda} 
		&= 
		[
			\begin{matrix}
				\alpha (z^{\lambda}_{1})^{\alpha -1}z^{\lambda}_{2} & (z^{\lambda}_{1})^{\alpha}
			\end{matrix}
		], \\
		Df^{\mu} 
		&= 
		[
			\begin{matrix}
				- (\alpha + 1)e^{\alpha \langle a,e_{1} \rangle + \langle a,e_{2} \rangle} (z^{\mu}_{1})^{\alpha} & \alpha (z^{\mu}_{2})^{\alpha-1}
			\end{matrix}
		], \\
		Df^{\sigma} 
		&= 
		[
			\begin{matrix}
				1 & -(\alpha +1)e^{\alpha \langle a,e_{1} \rangle + \langle a,e_{2} \rangle} (z^{\sigma}_{2})^{\alpha}
			\end{matrix}
		],
	\end{align*}
	respectively, we see the rank of $Df^{\lambda}$ and $Df^{\sigma}$ is one. However, the rank of $Df^{\mu}$ becomes zero at the point $(z^{\mu}_{1},z^{\mu}_{2}) = (0,0) \in \overline{C_{\mu}(V)}$. 
\end{example}

\begin{example}
	\label{eg: nonexampleA-1cp2}
	Let $X = \mathbb{C}P^{2}$, $k = 1$, and $V = \mathbb{R}p + a \; (a \in \mathbb{R}^{2})$ be an affine subspace spanned by $p = {^t[\alpha ~ -1]}$ for all integers $\alpha$ greater than or equal to three.
	Then, we can choose a basis $q$ of the orthogonal subspace to $V$ as $q = {^t [1 ~ \alpha]}$. 
	In this case, $\overline{C(V)}$ is not a complex submanifold in $X$. 
	Indeed, we give $f^{\lambda},f^{\mu},f^{\sigma}$ by 
	\begin{align*}
		f^{\lambda} 
		&= 
		z^{\lambda}_{1} (z^{\lambda}_{2})^{\alpha} - e^{\langle a,e_{1} \rangle + \alpha \langle a,e_{2} \rangle}, \\
		f^{\mu} 
		&= 
		z^{\mu}_{2} - e^{\langle a,e_{1} \rangle + \alpha \langle a,e_{2} \rangle} (z^{\mu}_{1})^{\alpha + 1}, \\
		f^{\sigma} 
		&= 
		(z^{\sigma}_{1})^{\alpha} - e^{\langle a,e_{1} \rangle + \alpha \langle a,e_{2} \rangle} (z^{\sigma}_{2})^{\alpha +1},
	\end{align*}
	respectively. 
	Note that $(0,0) \notin \overline{C_{\lambda}(V)}$.
	Since the Jacobian matrices are expressed as
	\begin{align*}
		Df^{\lambda} 
		&= 
		[
			\begin{matrix}
				(z^{\lambda}_{2})^{\alpha} & \alpha z^{\lambda}_{1} (z^{\lambda}_{2})^{\alpha -1}
			\end{matrix}
		], \\
		Df^{\mu} 
		&= 
		[
			\begin{matrix}
				- (\alpha + 1)e^{\langle a,e_{1} \rangle + \alpha \langle a,e_{2} \rangle} (z^{\mu}_{1})^{\alpha} & 1
			\end{matrix}
		], \\
		Df^{\sigma} 
		&= 
		[
			\begin{matrix}
				\alpha (z^{\sigma}_{1})^{\alpha -1} & - (\alpha + 1) e^{\langle a,e_{1} \rangle + \alpha \langle a,e_{2} \rangle} (z^{\sigma}_{2})^{\alpha}
			\end{matrix}
		],
	\end{align*}
	respectively, we see the rank of $Df^{\lambda}$ and $Df^{\sigma}$ is one. However, the rank of $Df^{\sigma}$ becomes zero at the point $(z^{\sigma}_{1},z^{\sigma}_{2}) = (0,0) \in \overline{C_{\sigma}(V)}$. 
\end{example}

By similar calculation, we see that $\overline{C(V)}$ is not a complex submanifold in $X$ if the slope of $V$ is not the same as treated above.

\begin{remark}
We can classify all examples for complex submanifolds $\overline{C(V)}$ in $X$ in terms of the conditions of $V$ by direct calculation.
In particular, when $X = \mathbb{C}P^{2}$, we can show that $\overline{C(V)}$ is a complex submanifold in $\mathbb{C}P^{2}$ if and only if 
the linear part of $V$ is spanned by ${^t[1 ~ 0]}$, ${^t[0 ~ 1]}$, ${^t[1 ~ 1]}$, ${^t[1 ~ 2]}$, ${^t[2 ~ 1]}$, or ${^t[1 ~ -1]}$.
\end{remark}

When $X = \mathbb{C}P^{2}$, we can determine the conditions that $\overline{C(V)}$ is a one-dimensional complex submanifold in $X$ by the linear part of an affine subspace $V = \mathbb{R}p+a$ in $\mathbb{R}^{2}$.

\subsection{Other Examples of Torus-equivariantly Embedded Toric Manifolds}
\label{subsec: other examples}

We demonstrate other examples of torus-equivariantly embedded toric manifolds in toric manifolds other than $\mathbb{C}P^{2}$.

\begin{figure}[hbtp]
	\centering
	\begin{tabular}{cc}
	\begin{minipage}[t]{0.5\hsize}
		\centering
		\includegraphics[keepaspectratio,width=5cm]{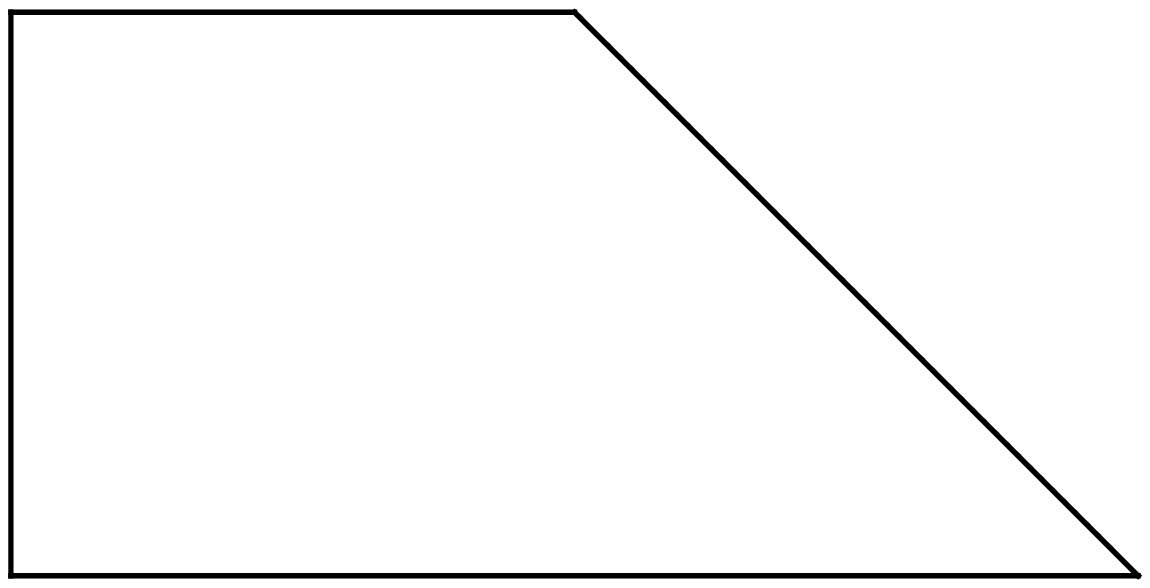}
		\caption{a Delzant polytope of $\mathbb{F}_{1}$}
		\label{fig: delzant f1}
	\end{minipage} &
	\begin{minipage}[t]{0.5\hsize}
		\centering
		\includegraphics[keepaspectratio,width=5cm]{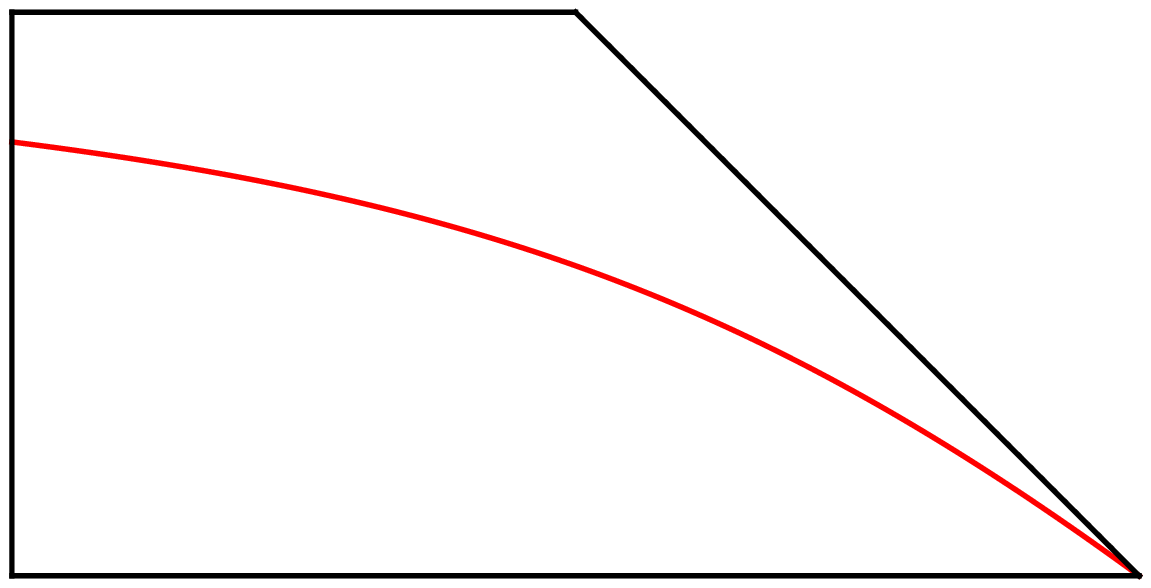}
		\caption{$\mu(\overline{C(V)})$ in Example \ref{eg: hirzebruch interior is complex submanifold}}
		\label{fig: example10f1}
	\end{minipage}
	\end{tabular}
\end{figure}

It is well-known that Delzant polytopes of $\mathbb{F}_{1}$ are shown in Figure \ref{fig: delzant f1}.
Define the points $\lambda$, $\mu$, $\sigma$, $\delta$ in $(\mathfrak{t}^{2})^{\ast} \cong \mathbb{R}^{2}$ by 
	\begin{equation*}
		\lambda = (0,0), ~
		\mu = (2,0), ~ 
		\sigma = (1,1), ~
		\delta = (0,1).
	\end{equation*}
Let $\Delta$ be a polytope defined by the convex hull of the points $\lambda$, $\mu$, $\sigma$, $\delta$.
Let $\Lambda = \{ \lambda, \mu,\sigma, \delta \}$ be a set of the vertices in the Delzant polytope of $\mathbb{F}_{1}$.
We define the inward pointing normal vectors to the facets by 
\begin{align*}
	&u^{\lambda}_{1} = 
		\left[
			\begin{matrix}
				1 \\ 0
			\end{matrix}
		\right],
	u^{\lambda}_{2} = 
		\left[
			\begin{matrix}
				0 \\ 1 
			\end{matrix}
		\right],
	u^{\mu}_{1} = 
		\left[
			\begin{matrix}
				0 \\ 1
			\end{matrix}
		\right],
	u^{\mu}_{2} = 
		\left[
			\begin{matrix}
				-1 \\ -1
			\end{matrix}
		\right], \\
	&u^{\sigma}_{1} = 
		\left[
			\begin{matrix}
				-1 \\ -1
			\end{matrix}
		\right],
	u^{\sigma}_{2} = 
		\left[
			\begin{matrix}
				0 \\ -1 
			\end{matrix}
		\right],
	u^{\delta}_{1} = 
		\left[
			\begin{matrix}
				0 \\ -1
			\end{matrix}
		\right],
	u^{\delta}_{2} = 
		\left[
			\begin{matrix}
				1 \\ 0
			\end{matrix}
		\right].
	\end{align*}
\begin{example}
	\label{eg: hirzebruch interior is complex submanifold}
	Let $X$ be a Hirzebruch surface $\mathbb{F}_{1}$ of degree one, $k=1$, and $p = {^t [1 ~ 0]}$.
	Then, we can choose a basis of the orthogonal subspace to $V$ as $q= {^t [0 ~ 1]}$.
	In this case, $\overline{C(V)}$ is a complex submanifold in $X$.
	Indeed, we give $f^{\lambda}$, $f^{\mu}$, $f^{\sigma}$, $f^{\delta}$ by 
	\begin{equation*}
		f^{\lambda}= z^{\lambda}_{1} - 1, \;
		f^{\mu} = z^{\mu}_{1} - z^{\mu}_{2}, \;
		f^{\sigma} = 1 - z^{\sigma}_{1}z^{\sigma}_{2}, \;
		f^{\delta} = 1 - z^{\delta}_{1},
	\end{equation*}
	respectively. 
	Note that $(0,0) \notin \overline{C_{\sigma}(V)}$.
	Since the Jacobian matrices are expressed as 
	\begin{equation*}
		Df^{\lambda} = 
		[\begin{matrix}
			0 & 1
		\end{matrix}], \;
		Df^{\mu} = 
		[\begin{matrix}
			1 & -1
		\end{matrix}], \;
		Df^{\sigma} = 
		[\begin{matrix}
			-z^{\sigma}_{2} & -z^{\sigma}_{1}
		\end{matrix}], \;
		Df^{\delta} = 
		[\begin{matrix}
			-1 & 0
		\end{matrix}],
	\end{equation*}
	respectively, we see the rank of each matrix is one.

	The image of $\mu\mid_{\overline{C(V)}}: \overline{C(V)} \to \mathbb{R}^{2}$ is given in Figure \ref{fig: example10f1}.
\end{example}

Delzant polytopes of $\mathbb{C}P^{3}$ are the convex hull of the points:
\begin{equation*}
	\lambda = (0,0,0), ~
	\mu = (2,0,0), ~
	\sigma = (0,0,2), ~ 
	\delta = (0,2,0).
\end{equation*}
We can see a Delzant polytope of a blow up of $\mathbb{C}P^{3}$ at the point corresponding to $\delta$ as the convex hull of the points:
\begin{align*}
	\lambda = (0,0,0), ~
	\mu = (2,0,0), ~ 
	\sigma = (0,0,2), \\
	\delta_{1} = (1,1,0), ~ 
	\delta_{2} = (0,1,0), ~
	\delta_{3} = (0,1,1).
\end{align*}
We define the inward pointing normal vectors to the facets by 
\begin{align*}
	&u^{\lambda}_{1} =
	\left[
		\begin{matrix}
			1 \\ 0 \\ 0
		\end{matrix}
	\right],
	u^{\lambda}_{2} =
	\left[
		\begin{matrix}
			0 \\ 1 \\ 0
		\end{matrix}
	\right],
	u^{\lambda}_{3} =
	\left[
		\begin{matrix}
			0 \\ 0 \\ 1
		\end{matrix}
	\right], 
	u^{\mu}_{1} = 
	\left[
		\begin{matrix}
			0 \\ 1 \\ 0
		\end{matrix}
	\right],
	u^{\mu}_{2} = 
	\left[
		\begin{matrix}
			-1 \\ -1 \\ -1
		\end{matrix}
	\right],
	u^{\mu}_{3} = 
	\left[
		\begin{matrix}
			0 \\ 0 \\ 1
		\end{matrix}
	\right], \\
	&u^{\sigma}_{1} = 
	\left[
		\begin{matrix}
			0 \\ 1 \\ 0
		\end{matrix}
	\right],
	u^{\sigma}_{2} = 
	\left[
		\begin{matrix}
			1 \\ 0 \\ 0
		\end{matrix}
	\right],
	u^{\sigma}_{3} = 
	\left[
		\begin{matrix}
			-1 \\ -1 \\ -1
		\end{matrix}
	\right], 
	u^{\delta_{1}}_{1} = 
	\left[
		\begin{matrix}
			0 \\ -1 \\ 0
		\end{matrix}
	\right],
	u^{\delta_{1}}_{2} = 
	\left[
		\begin{matrix}
			0 \\ 0 \\ 1
		\end{matrix}
	\right],
	u^{\delta_{1}}_{3} = 
	\left[
		\begin{matrix}
			-1 \\ -1 \\ -1
		\end{matrix}
	\right], \\
	&u^{\delta_{2}}_{1} = 
	\left[
		\begin{matrix}
			1 \\ 0 \\ 0
		\end{matrix}
	\right],
	u^{\delta_{2}}_{2} = 
	\left[
		\begin{matrix}
			0 \\ 0 \\ 1
		\end{matrix}
	\right],
	u^{\delta_{2}}_{3} = 
	\left[
		\begin{matrix}
			0 \\ -1 \\ 0
		\end{matrix}
	\right], 
	u^{\delta_{3}}_{1} = 
	\left[ 
		\begin{matrix}
			1 \\ 0 \\ 0
		\end{matrix}
	\right],
	u^{\delta_{3}}_{2} = 
	\left[
		\begin{matrix}
			0 \\ -1 \\ 0
		\end{matrix}
	\right],
	u^{\delta_{3}}_{3} = 
	\left[
		\begin{matrix}
			-1 \\ -1 \\-1
		\end{matrix}
	\right].
\end{align*}

\begin{example}
	\label{eg: blow up of cp3}
	Let $X$ be a blow up of $\mathbb{C}P^{3}$ at the point corresponding to $\delta$, $k = 2$, $p_{1} = {^t[1 ~ 0 ~ -1]}$ and $p_{2} = {^t[0 ~ 1 ~ 0]}$.
	Then, we can choose a basis of the orthogonal subspace to $V$ as $q = {^t [1 ~ 0 ~ 1]}$.
	In this case, $\overline{C(V)}$ is a complex submanifold in $X$.
	Indeed, we give $f^{\lambda}$, $f^{\mu}$, $f^{\sigma}$, $f^{\delta_{1}}$, $f^{\delta_{2}}$, $f^{\delta_{3}}$ by 
	\begin{align*}
		&f^{\lambda} = z^{\lambda}_{1}z^{\lambda}_{3} -1,
		f^{\mu} = z^{\mu}_{3} - (z^{\mu}_{2})^{2}, 
		f^{\sigma} = z^{\sigma}_{2} - (z^{\sigma}_{3})^{2}, \\
		&f^{\delta_{1}} = z^{\delta_{1}}_{2} - (z^{\delta_{1}}_{3})^{2},
		f^{\delta_{2}}= z^{\delta_{2}}_{1} z^{\delta_{2}}_{2} - 1,
		f^{\delta_{3}} = z^{\delta_{3}}_{1} - (z^{\delta_{3}}_{3})^{2},
	\end{align*}
	respectively.
	Note that $(0,0,0) \notin \overline{C_{\lambda}(V)}$, and $(0,0,0) \notin \overline{C_{\delta_{2}}(V)}$.
	Since the Jacobian matrices are expressed as 
	\begin{align*}
		&Df^{\lambda} = [z^{\lambda}_{3} ~ 0 ~ z^{\lambda}_{1}],
		Df^{\mu} = [0 ~ -2 z^{\mu}_{2} ~ 1], 
		Df^{\sigma} = [0 ~ 1 ~ -2z^{\sigma}_{3}], \\ 
		&Df^{\delta_{1}} = [0 ~ 1 ~ -2z^{\delta_{1}}_{3}], 
		Df^{\delta_{2}} = [z^{\delta_{2}}_{2} ~ z^{\delta_{2}}_{1} ~ 0], 
		Df^{\delta_{3}} = [1 ~ 0 ~ -2 z^{\delta_{3}}_{3}],
	\end{align*}
	respectively, we see the rank of each matrix is one.
\end{example}

In Example \ref{eg: blow up of cp3}, 
we have to suppose that $X$ is a blow up of $\mathbb{C}P^{3}$ at the point corresponding to the vertex $\delta$ because $\overline{C(V)}$ is not a complex submanifold in $X = \mathbb{C}P^{3}$ when the linear part of an affine subspace $V$ is spanned by $p_{1} = {^t[1 ~ 0 ~ -1]}$ and $p_{2} = {^t[0 ~ 1 ~ 0]}$.

\begin{ack}
The author is grateful to the advisor, Manabu Akaho for a lot of suggestions and supports. 
The author would also like to thank Yuichi Kuno and Yasuhito Nakajima for helpful discussions. 
This work is supported by JST, the establishment of university fellowships towards the creation of science technology innovation, Grant Number JPMJFS2139.
\end{ack}

\bibliography{equivariant_toric}
\bibliographystyle{alpha}

\end{document}